\documentclass[10pt]{article}

\usepackage[margin=1in]{geometry}
\usepackage{amsmath,amssymb,amsthm}
\theoremstyle{plain}
\newtheorem{theorem}{Theorem}
\newtheorem{lemma}[theorem]{Lemma}
\newtheorem{proposition}[theorem]{Proposition}
\newtheorem{corollary}[theorem]{Corollary}
\theoremstyle{definition}
\newtheorem{assumption}[theorem]{Assumption}
\newtheorem{remark}[theorem]{Remark}
\usepackage{booktabs}
\usepackage{graphicx}
\usepackage{caption}
\usepackage{microtype}
\usepackage[hidelinks]{hyperref}
\usepackage[section]{placeins}
\usepackage{authblk}

\title{\textbf{Optimal Stratified Allocation for Rare-Event Onset Forecasting in Dependent Sequences}}

\author[]{Jaskaran Singh}
\affil[]{Indraprastha Institute of Information Technology Delhi}
\date{September 2026}

\begin{document}
\maketitle

\begin{abstract}
\noindent
Let a finite population of $n$ labelled examples carry a class-weighted loss, with $N_1=\pi n$
examples in a rare positive class to which a learner applies the multiplier $\omega_1=N_0/N_1$. We
study estimation of the total risk from a subsample of size $K\ll n$ under designs that allocate
$K_0$ and $K_1$ draws to the two strata, and the transfer of that allocation to a model fitted on the
subsample. Three groups of results are obtained. First, we give the exact finite-population variance
of the weighted risk estimator under class-conditional sampling without replacement and solve for the
Neyman-optimal allocation $K_c^\star\propto N_cS_c$. The multiplier inflates the positive stratum's
dispersion by exactly the imbalance ratio, so the optimal allocation ratio is itself free of that
ratio and equal allocation displaces proportional allocation as the correct default; where the
classes are equidispersed before weighting, its efficiency exceeds proportional allocation by
$1/\{4\pi(1-\pi)\}$. Simple random sampling is dominated by an explicit between-stratum term, and an
exact bias identity shows that deterministic cluster-representative selection admits no unbiasedness
statement. A Serfling bound transfers the allocation to selection error over a finite candidate set.
Under the truncation $K_c=\min\{N_c,\lfloor K/2\rfloor\}$ the realised allocation ratio is
$\gamma=\min\{2\pi/f,1\}$, free of $n$, which yields a parameter-free predicted efficiency
$A(\pi,f)=\gamma/\{\pi(1-\pi)(1+\gamma)^2\}$ at every design point. Second, for sequentially indexed
data a measurability lemma bounds the record a labelled example occupies. That bound converts the
separation gap between training and evaluation blocks into a derived quantity, and the same quantity
limits the departure from block independence to $\beta(g-w-h-\Delta)$ under absolute regularity.
Third, ranking a score within the empirical distribution of the same model's held-out scores
preserves every within-block ordering, so its effect on pooled performance across separately fitted
models measures only the misalignment between them; with the $(n_V+1)$ correction the map becomes a
conformal $p$-value, valid up to $2\beta(\cdot)$.

We test all three against a detection problem that supplies the constants. The Phillips--Shi--Yu
procedure marks, after the fact, the dates on which a price series is too accelerated to be
consistent with a random walk, and central banks publish exuberance readings built from detectors of
this family. Our forward version of that question ranks securities by the probability that the
detector begins firing within $h\in\{5,10,20\}$ trading days, on a panel of 350 U.S.\ equities over
2004--2011 with fewer than 1\% positive rows and five purged forward evaluation blocks. Average
precision reaches 5.8, 3.1, and 2.3 times prevalence; the rank map quadruples pooled cross-block
average precision; and at a validation-fixed operating point the system attains 42.5\% episode recall
with a 6.5-day median lead at 1.64 false alerts per security-year. The predicted ordering of the four
designs holds, and at $h=10$ the ordering of the five design points predicted by $A(\pi,f)$ is
reproduced exactly (Spearman $\rho_S=1$, exact $p=0.0167$). The predicted dependence on $\pi$ across
horizons is not, and we identify the channels that lie outside the design-based argument instead of
absorbing the discrepancy into it. We run no backtest and make no claim of profitability.
\end{abstract}

\vspace{4pt}

\section{Introduction}

Two estimation problems recur whenever a learner is fitted repeatedly to a long, dependent sequence
in which the event of interest is rare. The first is allocation: given a budget of $K$ rows out of
$n$, how should the draws be divided between an abundant negative class and a scarce positive one,
and does the answer change when the loss already reweights the classes? The second is comparability:
scores produced by models fitted at different times lie on different scales, so a ranking that is
sound within one evaluation block need not survive pooling across blocks. Both admit exact answers
under a finite-population sampling model, and Section~\ref{sec:theory} gives them. Both are
questions about estimator variance and about measurability, not about the domain the data come from,
which is why the statements hold for any sequentially indexed panel with a rare binary target.

The instance that motivates them is a detection problem, and it fixes the constants that make the
theory testable. A price series occasionally stops behaving like ordinary fluctuation and begins to
accelerate, each increase raising the probability of the next. The Phillips--Shi--Yu (PSY) procedure
\cite{psy2015a,psy2015b} identifies such stretches after the fact, scanning a series day by day and
reporting whether the recent trajectory is too accelerated to be consistent with a random walk; its
output is a sequence of flagged and not-flagged dates. Section~\ref{sec:background} states the
procedure from first principles.

\paragraph{Why anticipating a detector's output is worth asking.}
Forecasting a statistical test may look like a strange target when one could forecast returns
directly. The answer is that this family of tests is already deployed as infrastructure. The
Federal Reserve Bank of Dallas, with Lancaster University and the University of Alicante,
computes exuberance indicators from these recursive right-tailed unit-root statistics and has
published them alongside every quarterly release of its International House Price Database since
the second quarter of 2013, covering more than two dozen countries
\cite{pavlidis2016,iho,dallasfed}; open implementations exist \cite{exuber}. For an institution
already acting on such a monitor, ``which series is about to start firing?'' is the operational
question rather than a proxy for a trading signal, because a monitor that reports only current
state offers no lead time for the response it is meant to inform. The target also has a useful
property: the underlying statistic is a published procedure rather than one we designed. Turning
it into a supervised label still requires four choices that we make, namely the 95\% critical
boundary, zero augmentation lags, a two-day merge gap, and a five-day minimum episode duration,
so the label is externally defined at its core, though not free of discretion at its edges.

By design, the PSY test is retrospective: it can say that a stretch of prices was accelerating
only after enough of that stretch has been observed. This paper asks the forward question. Let
$h \in \{5,10,20\}$ be a forecast horizon in trading days. At the close of day $t$, for every
stock the test is not currently flagging, we rank stocks by the probability that the test will
start flagging them somewhere in days $t+1$ through $t+h$. The ranking uses price and volume
information available by the close of day $t$, and no output of the test itself, since supplying
the test's own state would let a classifier reconstruct the detector rather than learn anything
antecedent to it.

\paragraph{Why the learning problem is hard.}
The task is hard for two reasons, and the theory addresses both. First,
information leaks backwards in time: a 60-day rolling feature is computed from a window of past
observations, so training and test examples close together in the index can share observations and a
model be scored on data it has partly seen. Lemma~\ref{lem:span} bounds that span and yields the
required separation, here 140 trading days against a bound of 87. Second, the positive class is
extremely rare, between 0.4\% and 1.4\% of eligible observations depending on horizon, so the
standard practice of training on everything and reweighting the rare class leaves an overwhelming
majority of near-duplicate negative rows. Whether a deliberately constructed smaller training set
helps, hurts, or merely saves time is exactly the allocation question above.

\medskip
\noindent\textbf{RQ1.} How well can information available at the close of day $t$ rank the
near-term risk that a stock enters a new flagged episode, under strict temporal separation?

\noindent\textbf{RQ2.} When the positive class is this rare, what is the optimal way to allocate a
fixed subset budget across the two classes, and does the allocation that theory selects improve the
accuracy--computation trade-off relative to training on everything?

\noindent RQ2 is the question the paper answers theoretically. Once the four constructions are seen
as sampling designs for the class-weighted empirical risk rather than as heuristics, the choice
between them is an allocation problem with an exact solution, and the empirical grid becomes a test
of a prediction rather than a search.

\paragraph{Contributions.}
\begin{enumerate}\itemsep2pt
\item An allocation theory for rare-event subset construction (Section~\ref{sec:theory}). We
show that the constructions compared here are stratified sampling designs for the class-weighted
empirical risk, give the exact finite-population variance of each, and solve for the Neyman-optimal
allocation. The default positive-class multiplier moves that optimum to equal class counts, with a
relative efficiency of $1/\{4\pi(1-\pi)\}$ over proportional allocation; uniform sampling is
dominated by an explicit between-class term; and clustering-based representative selection is placed
outside the framework by an exact bias identity rather than by measurement. A Serfling bound carries
the allocation to split-selection error at a node.
\item Two structural results about the evaluation. A measurability lemma bounds the record a
labelled example occupies, turning the purge length into a derived quantity that also bounds the
departure from block independence under a mixing assumption. A rank invariance result shows the
validation-relative score map cannot change any within-period ranking metric, so the entire pooled
cross-period gain is a measurement of misalignment between retraining dates.
\item A leakage-controlled forecasting task built on an established detector, in which examples
already inside a flagged episode are excluded, the label window begins strictly after the feature
date, and all blocks are separated by a purge satisfying the bound above
(Sections~\ref{sec:problem} and \ref{sec:protocol}).
\item An operational early-warning evaluation converting rankings into alerts under a false-alarm
budget fixed on validation data, reporting detection rate and warning lead jointly with the alert
burden they cost (Section~\ref{sec:alerts}).
\item A validation-relative percentile score making rankings from models retrained at different
dates comparable without test labels, which roughly quadruples pooled cross-period average precision
(Section~\ref{sec:percentile}).
\item A documented, monitorable failure mode: in one test period the ranking is inverted, a
covariate-shift monitor does not flag it, a label-conditional sign reversal in a specific feature
cluster does, and we propose a validation-time gating rule on that basis
(Section~\ref{sec:inversion}).
\item An empirical test of the allocation theory across 54 configurations. The predicted ordering of
the four rules holds: class-balanced sampling improves average precision over full-data training at
every 10-day retention level, while uniform, class-stratified, and $K$-means construction show no
comparable improvement. The predicted dependence of the margin on rarity does not hold, and we report
that failure and the three channels outside the design-based argument that could account for it
(Sections~\ref{sec:th-predict}, \ref{sec:reduction} and \ref{sec:reduction-results}).
\end{enumerate}

\paragraph{Scope.}
The event forecast is the output of a statistical test, not an economic fact: the test can fire on a
stretch of no lasting consequence, and a real speculative episode can fail to trigger it. We run no
backtest and make no claim of profitability.

\section{Background: when does a price count as explosive?}
\label{sec:background}

The null model for an unremarkable log price $p_t$ is the random walk
$p_t = p_{t-1} + \varepsilon_t$, with $\varepsilon_t$ mean-zero: past movement carries no information
about future direction. The alternative is $p_t = (1+\delta)p_{t-1} + \varepsilon_t$ with $\delta>0$,
under which deviations compound; $\delta = 0$ is the unit-root boundary case. The Augmented
Dickey--Fuller (ADF) test discriminates the two on a fixed window by regressing
$r^{(1)}_\tau = p_\tau - p_{\tau-1}$ on $p_{\tau-1}$ and reporting the $t$-statistic on that
coefficient, right-tailed, so only large positive values count as evidence of acceleration.

A single window cannot locate an episode's start. The PSY procedure runs the ADF regression on every
window ending at $t$ and beginning at each admissible earlier date, retains the largest statistic,
and compares it against a threshold simulated separately for each window length, since the null
distribution depends on the amount of data available. This is the backward supremum ADF (BSADF)
statistic, and it yields a flag on every trading day. The procedure reads one price series and has
access to no news, valuation, or order-flow information, so its output is a hypothesis-test rejection
rather than a diagnosis. We use it as a reproducible label generator with an established literature
behind it, not as a claim about economic reality.

\section{Preliminaries and notation}
\label{sec:notation}

Securities are indexed by $i$ and rows by $t$, a position in that security's own history rather than
a calendar date. $\mathbb{I}\{A\}$ is the indicator, $\Pr(\cdot)$ probability, $\mathbb{E}[\cdot]$
expectation. A hat marks a quantity computed on a test block, a superscript $(h)$ dependence on the
horizon, and a superscript $V$ a quantity from a validation block. Each symbol below carries exactly
one meaning throughout. Pure indices are bound by the expression containing them: $j$ ranges over rows of whatever set is
named, and $\theta$ over score thresholds; $k$ is a lag length in Table~\ref{tab:features}, a feature
index in \eqref{eq:reversal}, and a day count in Section~\ref{sec:limitations}. Where a
statement is about distributions rather than about a particular row, $X$ and $Y$ denote a generic
feature vector and its label, as in $\Pr(Y\mid X)$. The superscript $\mathrm{full}$ marks a quantity
from the full-data run of the same cell, $L_2$ is the squared-norm penalty on leaf weights, and
$z$-score and $z$-statistic carry their usual meanings, namely a value standardised by a mean and a
standard deviation.

\begin{table}[htbp]
\centering
\caption{Symbols. Panel A: data, detector, and learning task. Panel B: sampling design, evaluation,
and alerts; Section~\ref{sec:theory} uses Panel B throughout.}
\label{tab:notation}
\small
\begin{tabular}{@{}l@{\hspace{1em}}p{0.70\textwidth}@{}}
\toprule
\multicolumn{2}{@{}l}{\textbf{A. Data, detector, learning task}}\\
\midrule
$i$, $t$, $\tau$ & Security; row within that security; time index inside one ADF window\\
$P_{i,t}$, $p_{i,t}=\log P_{i,t}$, $V_{i,t}$ & Adjusted close; log price; volume ($\bar V$, $s_V$ its rolling mean and s.d.)\\
$r^{(1)}_{i,t}$, $\varepsilon_t$, $\delta$ & One-day log return; mean-zero shock; explosive coefficient in $p_t=(1+\delta)p_{t-1}+\varepsilon_t$\\
$\mathrm{BSADF}_{i,t}$, $\mathrm{cv}^{0.95}_{i,t}$, $m_i$ & Backward supremum ADF statistic; simulated threshold; minimum window\\
$R_{i,t}$, $E_{i,t}$, $O_{i,t}$, $W$ & Raw flag; retained-episode membership; onset indicator; waiting time to next onset\\
$h$, $Y^{(h)}_{i,t}$, $T_i$ & Horizon $h\in\{5,10,20\}$; label \eqref{eq:target}; rows available for security $i$\\
$x_{i,t}$, $x_{i,t,k}$, $d$ & Feature vector; its $k$th entry; feature columns (17, of which 16 distinct)\\
$q$, $q'$, $u(q)$, $\tilde u(q)$ & Score; prevalence-shifted score; validation percentile \eqref{eq:percentile}; corrected form \eqref{eq:conformal}\\
$\alpha$, $\beta$, $q^V_j$, $n_V$ & Calibration intercept and slope; $j$th validation score; their count ($\bar q^V$, $\bar y^V$ their mean and event rate)\\
$\Phi_M$, $\phi_m$, $M$, $\eta$, $\lambda_\kappa$ & Ensemble and $m$th tree; rounds; learning rate; hazard (Section~\ref{sec:limitations} only)\\
\addlinespace
\multicolumn{2}{@{}l}{\textbf{B. Sampling design, evaluation, alerts}}\\
\midrule
$n$, $N_0$, $N_1$, $\pi$ & Rows in the training fold; negative and positive counts; event rate $N_1/n$\\
$f$, $K$, $K^*$, $K_c$, $S_c$ & Retention fraction; size $\lceil fn\rceil$; floored $\max(K,50)$; rows kept from class $c$; that set\\
$w_j$, $\ell_j$, $\omega_c$ & Weight of row $j$; its class-weighted loss \eqref{eq:lossrow}; multiplier, $\omega_1=N_0/N_1$\\
$S_c^2$, $\sigma_c^2$, $B$ & Within-class dispersion of $\ell$ \eqref{eq:strvar} and of the unweighted loss; between-class term\\
$a_c$, $a^\star_c$, $\gamma$ & Realised and Neyman-optimal budget shares; realised allocation ratio \eqref{eq:gamma}\\
$w$, $\delta_{\mathrm{m}}$, $d_{\min}$, $\Delta$, $g$ & Feature window; merge gap; minimum duration; $\Delta=\delta_{\mathrm{m}}+d_{\min}$; purge length\\
$\beta(\cdot)$, $I$, $B_{\mathrm{mb}}$ & Absolute regularity coefficients; MiniBatch iterations and batch size\\
$F$, $F_0,\dots,F_4$, $s$ & A fold; the five forward folds in order; seed, $s\in\{11,23,47\}$\\
$\mathrm{AP}$, ROC-AUC, $\mathrm{prec}_\theta$, $\mathrm{rec}_\theta$ & Average precision; area under ROC; precision and recall at the $\theta$th threshold\\
$\Delta^{(h,f)}_{F,s}$, $\hat\pi_+$, $\rho^{(h)}_F(k)$ & Paired AP difference \eqref{eq:paired}; fraction of draws positive; feature--label correlation \eqref{eq:reversal}\\
$b$, $t_o$, $C_{\text{false}}$, $N_{\text{elig}}$ & Alert budget per security-year; onset date; false clusters; eligible rows\\
\bottomrule
\end{tabular}
\end{table}

\section{Related work}

\paragraph{Date-stamping explosive behaviour, and its use as live infrastructure.}
Phillips et al.\ \cite{pwy2011} introduced recursive right-tailed testing for locating episodes of
exuberance, and \cite{psy2015a,psy2015b} generalised it to multiple episodes within one series. The
statistics are used operationally: \cite{pavlidis2016} applies them to the Dallas Fed International
House Price Database, the indicators are released quarterly \cite{dallasfed} and published as a
cross-country monitor \cite{iho}, and open implementations exist \cite{exuber}. We take the detected
state as given and add a forecast made strictly before onset. Prior work pairing explosive-root
indicators with machine learning covers the S\&P 500 and metal prices \cite{basoglu2021,ozgur2021};
relative to that line we add a prospective exclusion rule, paired rare-event subset comparisons with
an allocation theory behind them, and an alert-budget analysis. Biagini et al.\ \cite{biagini2025}
use option prices and deep learning under a local-martingale definition of a bubble, a different
target and framing, so we cite it as related design rather than as a baseline. Harvey et al.\
\cite{harvey2016} show this test family is non-pivotal and over-sized under permanent volatility
shifts, which Section~\ref{sec:limitations} treats as a specific threat to our labels.

\paragraph{Sampling, rare-event evaluation, and temporal validity.}
Section~\ref{sec:theory} is finite-population sampling theory in the tradition of Neyman allocation
for stratified designs, applied to the empirical risk of a class-weighted learner; the concentration
step uses Serfling's inequality for sampling without replacement \cite{serfling1974}, and the
class-balanced design is a case-control scheme, equivalent to random undersampling with importance
weights. Bachem et al.\ \cite{bachem2018} give approximation guarantees for clustering-shaped
objectives, which Proposition~\ref{prop:kmeans} shows do not transfer to a loss that varies with the
model state. When positives are rare, precision--recall summaries are more informative than ROC-AUC
\cite{saito2015}; LightGBM \cite{lightgbm} is our nonlinear benchmark. We do not evaluate
gradient-based one-side sampling, which reduces data inside the boosting loop rather than before it,
and Section~\ref{sec:limitations} flags the omission. Random cross-validation is inappropriate under
temporal dependence; L\'opez de Prado \cite{lopezdeprado2018} formalises purging, which
Lemma~\ref{lem:span} recovers as a measurability condition, and calibration is fitted on the
validation block by Platt's method \cite{platt1999}. Equation~\eqref{eq:percentile} is structurally
the split-conformal transformation \cite{vovk2005} applied per retraining window; the dependent-case
bound of Proposition~\ref{prop:conformal} follows the coupling route of \cite{doukhan1994}. For the
anomalous period of Section~\ref{sec:inversion} we use the shift vocabulary of \cite{morenotorres2012}.

\section{Problem formulation}
\label{sec:problem}

Let $i$ index securities and $t$ index trading rows within a security. Write $P_{i,t}$ for the
split- and dividend-adjusted closing price, $p_{i,t} = \log P_{i,t}$, and $V_{i,t}$ for daily
volume.

\paragraph{Label construction.}
Let $\mathrm{BSADF}_{i,t}$ be the statistic of Section~\ref{sec:background} and $\mathrm{cv}^{0.95}_{i,t}$
its simulated 95th-percentile null threshold. The raw flag is
$R_{i,t} = \mathbb{I}\{\mathrm{BSADF}_{i,t} > \mathrm{cv}^{0.95}_{i,t}\}$. Raw flags flicker around the
boundary, so we apply two deterministic filters: runs separated by at most two non-flagged days
are merged, and a merged run is retained only if it lasts at least five days. Let $E_{i,t}=1$
indicate membership in a retained episode. The onset indicator is
\begin{equation}
O_{i,t} = \mathbb{I}\{E_{i,t}=1 \text{ and } E_{i,t-1}=0\},
\end{equation}
one on the first day of an episode and zero elsewhere.

\paragraph{Forecasting target.}
For horizon $h$,
\begin{equation}
Y^{(h)}_{i,t} = \mathbb{I}\Big\{\textstyle\sum_{j=1}^{h} O_{i,t+j} \ge 1\Big\}.
\label{eq:target}
\end{equation}
The sum starts at $j=1$, so every forecast onset lies strictly after the feature date. A row is
eligible only if
\begin{equation}
E_{i,t}=0,\qquad t+h \le T_i,\qquad x_{i,t}\ \text{fully observed},
\label{eq:eligible}
\end{equation}
where $T_i$ is the number of available rows for security $i$ and $x_{i,t}$ the feature vector. The
first condition removes rows inside an ongoing episode, so the model is never asked to recognise
something already under way. One episode contributes up to $h$ positive rows, so row-level metrics
overstate the number of independent events, which is why we also report episode-level results
(Section~\ref{sec:alerts}).

\paragraph{Eligibility uses a filtered indicator.}
Since $E_{i,t}$ depends on the merge and duration filters it is knowable only after the fact: a live
monitor knows $R_{i,t}$ but not $E_{i,t}$ and would score a slightly different set of securities. No
future information reaches the inputs, every feature still being built from data no later than $t$;
what differs is the evaluation population. Across the 705{,}061 rows the detector raises 27{,}578 raw
flags and retains 26{,}350, and summing per security the excess of raw over retained days bounds the
disagreement below by 1{,}573 rows, 0.22\% of the panel and 5.7\% of raw-flagged rows. Filtering both
removes flagged days in short runs and adds unflagged days inside merged gaps, so this does not fix
the direction of the effect. Defining eligibility on $R_{i,t}=0$ is observable at $t$ and is the first
change we would make. The second condition in \eqref{eq:eligible} censors incomplete horizons rather
than labelling them negative.

\section{Theory}
\label{sec:theory}

This section states the results the empirical sections test. Section~\ref{sec:th-span} bounds the
stretch of the record a single labelled example occupies, which fixes the purge length and supplies
the separation hypothesis the later statements use. Section~\ref{sec:th-design} recasts the four
training-subset constructions of Section~\ref{sec:reduction} as sampling designs for the
class-weighted empirical risk and gives the exact finite-population variance of each.
Section~\ref{sec:th-alloc} solves for the optimal allocation and shows what the class-weighting
multiplier does to it. Section~\ref{sec:th-predict} converts these into predictions, including one
the data does not support. Section~\ref{sec:th-align} treats the cross-period score map. Proofs are
in Appendix~\ref{app:proofs}.

\subsection{The information span of a labelled example}
\label{sec:th-span}

Let $\mathcal{G}_{[a,b]}$ denote the $\sigma$-field generated by all prices and volumes on market
dates $a$ through $b$. Write $w$ for the deepest feature window, $\delta_{\mathrm{m}}$ for the merge gap, and
$d_{\min}$ for the minimum episode duration, and set $\Delta = \delta_{\mathrm{m}} + d_{\min}$.

\begin{lemma}[Information span]
\label{lem:span}
For every $(i,t)$ and every $h$, the feature vector $x_{i,t}$ is $\mathcal{G}_{[t-w,\,t]}$-measurable
and the label $Y^{(h)}_{i,t}$ is $\mathcal{G}_{(-\infty,\,t+h+\Delta]}$-measurable. Consequently, if
the last training date is $t_{\mathrm{tr}}$, the first test date is $t_{\mathrm{te}}$, and
$g = t_{\mathrm{te}} - t_{\mathrm{tr}}$, then no test example reads any date that a training label
reads, provided
\begin{equation}
g \;>\; w + h + \Delta .
\label{eq:purgebound}
\end{equation}
\end{lemma}

The label is not measurable with respect to a bounded window on the left, because
$\mathrm{BSADF}_{i,t}$ regresses over windows reaching back to the start of the series. That
asymmetry is harmless: the backward reach of a training label points into data the training block
already owns. What matters is the forward reach of a training label against the backward reach of a
test feature, and \eqref{eq:purgebound} is exactly that comparison. With $w = 60$, $\delta_{\mathrm{m}} = 2$,
$d_{\min} = 5$, and $h \le 20$, the bound is $87$ trading days, and the design uses $g = 140$. This
replaces the informal count of Section~\ref{sec:protocol} with a statement about which dates enter
which $\sigma$-field.

Non-overlap is not independence, since the BSADF statistic depends on the whole past. Independence
requires a mixing condition, which we assume rather than establish.

\begin{assumption}[Mixing]
\label{ass:mixing}
The date-indexed panel process is absolutely regular with coefficients $\beta(\cdot)$ satisfying
$\beta(u) \to 0$ as $u \to \infty$.
\end{assumption}

\begin{lemma}[Approximate block independence]
\label{lem:mixing}
Under Assumption~\ref{ass:mixing} and \eqref{eq:purgebound}, for bounded measurable $\varphi$ of the
training block and $\psi$ of the test block,
\begin{equation}
\big|\,\mathbb{E}[\varphi\psi] - \mathbb{E}[\varphi]\,\mathbb{E}[\psi]\,\big|
\;\le\; \beta\big(g - w - h - \Delta\big)\,\|\varphi\|_\infty\|\psi\|_\infty .
\label{eq:mixgap}
\end{equation}
\end{lemma}

Every subsequent statement about a test block is exact under independence and degrades by the
right-hand side of \eqref{eq:mixgap} otherwise. We do not estimate $\beta$, and
Section~\ref{sec:limitations} records that this leaves the degradation term unquantified.

\subsection{Subset construction as a sampling design}
\label{sec:th-design}

Fix a fold and index its $n$ training rows by $i$, with $N_1$ positives and $N_0 = n - N_1$
negatives at the horizon in question, and event rate $\pi = N_1/n$. Fix a candidate model state
$\Phi$ and write
\begin{equation}
\ell_i \;=\; \omega_{Y_i}\,\ell\big(\Phi(x_i),\,Y_i\big),
\qquad \omega_1 = N_0/N_1,\quad \omega_0 = 1,
\label{eq:lossrow}
\end{equation}
for the class-weighted loss contribution of row $i$, with $\omega_1$ the multiplier of
Section~\ref{sec:models} computed from the complete fold. The full-data objective is
$L(\Phi) = \sum_{i=1}^{n}\ell_i$. A construction rule returns $S = S_0 \cup S_1$ with
$|S_c| = K_c$ and weights $w_j$, and the model is fitted to
$\widehat L(\Phi) = \sum_{j \in S} w_j \ell_j$. Write $\bar\ell_c$ for the class-$c$ mean of
$\ell$, $\bar\ell$ for the overall mean, and
\begin{equation}
S_c^2 \;=\; \frac{1}{N_c-1}\sum_{i:\,Y_i=c}\big(\ell_i - \bar\ell_c\big)^2
\label{eq:strvar}
\end{equation}
for the within-class dispersion of the weighted loss. Let $\sigma_c^2$ be the same quantity computed
from the unweighted contributions $\ell(\Phi(x_i),Y_i)$, so that $S_1 = \omega_1\sigma_1$ and
$S_0 = \sigma_0$.

Three of the four rules are probability samples: \emph{uniform} draws $K$ rows from all $n$ without
replacement with $w_j = n/K$; \emph{class-stratified} and \emph{class-balanced} draw $K_c$ rows
without replacement within class $c$ with $w_j = N_c/K_c$, differing only in how $K_c$ is set. The
fourth, $K$-means representatives, is not.

\begin{proposition}[Unbiasedness]
\label{prop:unbiased}
For each of the three probability designs and every fixed $\Phi$,
$\mathbb{E}\big[\widehat L(\Phi)\big] = L(\Phi)$.
\end{proposition}

This generalises Equation~\eqref{eq:unbiased}, which is the class-stratified case, and makes explicit
that unbiasedness holds pointwise in $\Phi$ and therefore says nothing on its own about which design
is preferable. The dispersion does.

\begin{proposition}[Exact finite-population variance]
\label{prop:variance}
For the class-conditional designs,
\begin{equation}
\operatorname{Var}\big(\widehat L(\Phi)\big)
\;=\; \sum_{c \in \{0,1\}} \frac{N_c^2}{K_c}\Big(1 - \frac{K_c}{N_c}\Big) S_c^2 ,
\label{eq:varstrat}
\end{equation}
and for the uniform design,
$\operatorname{Var}(\widehat L) = \frac{n^2}{K}\big(1 - \frac{K}{n}\big)S^2$ with $S^2$ the pooled
dispersion of $\ell$ over all $n$ rows.
\end{proposition}

\begin{proposition}[Uniform sampling is dominated]
\label{prop:unifdom}
Let $B = \sum_c N_c(\bar\ell_c - \bar\ell)^2$ denote the between-class term. Then, at a common
budget $K$ and with $K_c = KN_c/n$ in the stratified design,
\begin{equation}
\operatorname{Var}_{\mathrm{unif}} - \operatorname{Var}_{\mathrm{prop}}
\;=\; \frac{n^2\,(1 - K/n)}{K\,(n-1)}
\Big[\,B \;+\; \tfrac{1}{n}\textstyle\sum_c N_c S_c^2 \;-\; \sum_c S_c^2 \,\Big],
\label{eq:unifgap}
\end{equation}
which is nonnegative whenever $B \ge \sum_c S_c^2$. Under \eqref{eq:lossrow} the two class means
differ by a factor of order $\omega_1 = N_0/N_1$, so the condition holds by a wide margin in the
rare-event regime.
\end{proposition}

The dominance is not a statement about sample size but about allocation noise: uniform sampling is
proportional allocation in which $K_1$ is itself random, hypergeometric with mean $K\pi$, and the
law of total variance charges for that randomness. At $f = 0.01$ and $\pi = 0.0037$ the expected
positive count in a uniform draw from a fold of $n$ rows is $0.0037\,K$, and the probability of a
draw containing no positive at all is approximately $(1-\pi)^K$.

\begin{proposition}[Representative selection admits no unbiasedness statement]
\label{prop:kmeans}
Let clusters $C_1,\dots,C_{K_c}$ partition class $c$, let $r(C)$ be the retained representative of
cluster $C$, and let $w_{r(C)} = |C|$. Then
\begin{equation}
\widehat L - L \;=\; \sum_{c}\sum_{C} |C|\,\big(\ell_{r(C)} - \bar\ell_C\big),
\label{eq:kmbias}
\end{equation}
which vanishes for all $\Phi$ only if $\ell$ is constant within every cluster. Clusters are formed
in feature space without reference to $\ell$, and $\ell$ changes at every boosting round, so no
choice of clustering makes \eqref{eq:kmbias} vanish uniformly in $\Phi$.
\end{proposition}

Proposition~\ref{prop:kmeans} places $K$-means representatives outside the allocation theory
entirely rather than at a poor point within it, which is a sharper statement than the empirical
observation that they underperform. Their construction cost, analysed in
Section~\ref{sec:kmeanscost}, therefore buys no compensating guarantee.

\subsection{Optimal allocation, and what the class multiplier does to it}
\label{sec:th-alloc}

\begin{theorem}[Neyman allocation for the class-weighted risk]
\label{thm:neyman}
Subject to $K_0 + K_1 = K$ with $0 < K_c \le N_c$, the allocation minimising \eqref{eq:varstrat} is
\begin{equation}
K_c^{\star} \;=\; K\,\frac{N_c S_c}{N_0 S_0 + N_1 S_1},
\label{eq:neyman}
\end{equation}
and for any feasible allocation with shares $a_c = K_c/K$ and optimal shares
$a^\star_c = K^\star_c/K$, the efficiency loss is
\begin{equation}
\frac{\operatorname{Var}(a)}{\operatorname{Var}(a^\star)}
\;=\; \sum_{c}\frac{(a^\star_c)^2}{a_c} \;\ge\; 1 ,
\label{eq:alloceff}
\end{equation}
with equality only at $a = a^\star$.
\end{theorem}

\begin{corollary}[Equal allocation is optimal under equal unweighted dispersion]
\label{cor:balopt}
Equal allocation $K_0 = K_1$ is Neyman-optimal if and only if $N_0S_0 = N_1S_1$, which under
\eqref{eq:lossrow} with $\omega_1 = N_0/N_1$ is equivalent to
\begin{equation}
\sigma_1 \;=\; \sigma_0 .
\label{eq:balcond}
\end{equation}
Proportional allocation $K_c \propto N_c$ is Neyman-optimal if and only if $S_0 = S_1$, equivalently
$\sigma_1 = (N_1/N_0)\,\sigma_0 = \{\pi/(1-\pi)\}\,\sigma_0$. Equivalently, and without any condition
on the dispersions,
\begin{equation}
\frac{K_1^\star}{K_0^\star} \;=\; \frac{N_1S_1}{N_0S_0} \;=\; \frac{\sigma_1}{\sigma_0},
\label{eq:ratioform}
\end{equation}
so under the class-weighted loss the optimal split between strata does not depend on $\pi$ at all.
\end{corollary}

The multiplier $\omega_1 = N_0/N_1$ that rare-event learners apply by default inflates the positive
rows' contributions by exactly the imbalance ratio, and thereby moves the Neyman optimum from
proportional allocation to equal allocation. Substituting into \eqref{eq:neyman} gives the ratio
form $K_1^\star/K_0^\star = \sigma_1/\sigma_0$, in which the imbalance ratio has cancelled: once the
loss is class-weighted, the optimal split between strata depends on their unweighted dispersions
alone. Equal allocation therefore needs only condition \eqref{eq:balcond}, while proportional
allocation would need the positive class to be less dispersed than the negative class by a factor of
order $\pi$, implausible in a fold with $\pi < 0.015$. Class-balanced construction and the loss
multiplier are complements, not two applications of one correction.

\begin{corollary}[Relative efficiency]
\label{cor:releff}
Under \eqref{eq:balcond} and ignoring the finite-population corrections, at a common budget $K$,
\begin{equation}
\frac{\operatorname{Var}_{\mathrm{prop}}}{\operatorname{Var}_{\mathrm{bal}}}
\;=\; \frac{n^2}{4N_0N_1} \;=\; \frac{1}{4\pi(1-\pi)} .
\label{eq:releff}
\end{equation}
\end{corollary}

At the three mean fold prevalences of Table~\ref{tab:macro}, $\pi = 0.0037$, $0.0071$, and $0.0143$,
the factor is $67.8$, $35.5$, and $17.7$. Combining Propositions~\ref{prop:unifdom} and
\ref{prop:kmeans} with Corollaries~\ref{cor:balopt} and \ref{cor:releff} gives a complete predicted
ordering of the four rules, which Section~\ref{sec:reduction-results} and
Appendix~\ref{app:grids} test.

The variance in \eqref{eq:varstrat} is that of the objective at a fixed $\Phi$, and boosting selects
$\Phi$ adaptively using $\widehat L$ itself. The following bridges the two at the cost of a
finite-candidate assumption.

\begin{proposition}[Uniform control over a finite candidate set]
\label{prop:uniform}
Let $\mathcal{C}$ be a finite set of candidate splits considered at a node, with
$|\mathcal{C}| = M$, and suppose the per-row contributions are bounded, $|\ell_i| \le \Lambda$. Then
for the class-conditional designs, with probability at least $1 - \varepsilon$,
\begin{equation}
\sup_{\Phi \in \mathcal{C}} \big|\widehat L(\Phi) - L(\Phi)\big|
\;\le\; \Lambda \sum_{c} N_c \sqrt{\frac{2\,(1 - (K_c-1)/N_c)}{K_c}\,\log\frac{2M}{\varepsilon}} ,
\label{eq:serfling}
\end{equation}
by Serfling's inequality for sampling without replacement together with a union bound.
\end{proposition}

The rate in \eqref{eq:serfling} is $\mathcal{O}(K_c^{-1/2})$ per class with the same $N_c$ weighting
that drives \eqref{eq:neyman}, so allocation governs split-selection error at the same margin at
which it governs objective variance. We do not claim a bound over the boosting path, only at a node
with a fixed candidate set; Section~\ref{sec:limitations} records this.

\subsection{What the theory predicts, and where it should fail}
\label{sec:th-predict}

The implemented rule sets $K_c = \min\{N_c, \max(25,\lfloor K^*/2\rfloor)\}$, so when positives are
scarce the design retains all of them and the realised allocation falls short of equal. The resulting
loss is exactly computable and, usefully, free of every unknown quantity.

\begin{corollary}[Truncated allocation, in closed form]
\label{cor:trunc}
With $K^* = \lceil fn\rceil$ and $N_1 = \pi n$, and ignoring the floor at $25$, the realised
allocation ratio is
\begin{equation}
\gamma \;=\; K_1/K_0 \;=\; \min\{2\pi/f,\ 1\},
\label{eq:gamma}
\end{equation}
which depends on $\pi$ and $f$ alone and not on the fold size $n$. By \eqref{eq:alloceff} the
efficiency loss relative to equal allocation is $(1+\gamma)^2/(4\gamma)$, and combining with
Corollary~\ref{cor:releff} the predicted efficiency of the implemented design relative to
proportional allocation is
\begin{equation}
A(\pi,f) \;=\; \frac{\gamma}{\pi\,(1-\pi)\,(1+\gamma)^2}.
\label{eq:netadv}
\end{equation}
\end{corollary}

Equation~\eqref{eq:netadv} assigns a number to every cell of the experimental grid before any model
is fitted, with no free parameter and no appeal to the data (Table~\ref{tab:predict}). Three
predictions follow, and the evidence treats them differently.

\begin{description}\itemsep2pt
\item[P1 (ordering of the rules).] Class-balanced construction should dominate class-stratified,
which should dominate uniform by the between-class term of Proposition~\ref{prop:unifdom}, while
$K$-means representatives are governed by the bias identity \eqref{eq:kmbias} rather than by any
allocation argument.
\item[P2 (between horizons).] $A(\pi,f)$ is decreasing in $\pi$ at fixed $f$, so the advantage should
be largest at $h=5$ and smallest at $h=20$.
\item[P3 (within horizon).] $A(\pi,f)$ is non-increasing in $f$ at fixed $\pi$, strictly so once
$f > 2\pi$, so within each horizon the advantage should decline across the retention grid.
\end{description}

\begin{table}[htbp]
\centering
\caption{Predicted efficiency $A(\pi,f)$ of the implemented class-balanced design relative to
proportional allocation, from Corollary~\ref{cor:trunc}, at the mean fold prevalences of
Table~\ref{tab:macro}; the realised ratio $\gamma$ of \eqref{eq:gamma} is in parentheses. Every entry
is fixed before estimation and uses no fitted quantity.}
\label{tab:predict}
\small
\begin{tabular}{clrrrrr}
\toprule
$h$ & $\pi$ & $f=1\%$ & $2.5\%$ & $5\%$ & $10\%$ & $20\%$\\
\midrule
5  & 0.0037 & 66.3 (0.74) & 47.8 (0.30) & 30.5 (0.15) & 17.4 (0.07) & 9.3 (0.04)\\
10 & 0.0071 & 35.5 (1.00) & 32.8 (0.57) & 24.4 (0.28) & 15.4 (0.14) & 8.8 (0.07)\\
20 & 0.0143 & 17.7 (1.00) & 17.7 (1.00) & 16.4 (0.57) & 12.3 (0.29) & 7.8 (0.14)\\
\bottomrule
\end{tabular}
\end{table}

Section~\ref{sec:reduction-results} finds P1 and P3 supported and P2 contradicted. We state the
contradiction rather than absorb it. Three channels lie outside the design-based argument and are
candidates. The variance in \eqref{eq:varstrat} governs the objective, not the generalisation gap,
and at $h = 5$ full-data training already attains ROC-AUC above $0.85$ in four of five folds, leaving
little headroom for a variance improvement to convert into a ranking improvement. The negative class
consists largely of temporally adjacent near-duplicate rows, so the effective $N_0$ entering
\eqref{eq:releff} overstates the independent information available, and does so by an amount that
varies with horizon; replacing $N_0$ by an effective sample size would compress the columns of
Table~\ref{tab:predict} in the direction the data require. And roughly half the measured $h = 10$
advantage arises in the single period $F_2$, where full-data training produces an inverted ranking;
that is a concept-shift event, and no allocation argument speaks to it.
Section~\ref{sec:limitations} states the experiments that would separate these.

\subsection{Score alignment across retraining windows}
\label{sec:th-align}

The map $u$ of Equation~\eqref{eq:percentile} is the empirical distribution function of the model's
own validation scores evaluated at the test score, hence nondecreasing.

\begin{proposition}[Within-fold rank invariance]
\label{prop:rankinv}
Fix a fold and let $q_1,\dots,q_N$ be its test scores. If no two of them fall in the same gap of the
validation score sample, so that $u$ is injective on $\{q_1,\dots,q_N\}$, then the within-fold
average precision and ROC-AUC of $\{u(q_r)\}$ equal those of $\{q_r\}$ exactly. Without injectivity
the two differ only through the ties $u$ creates.
\end{proposition}

\begin{corollary}
\label{cor:alignment}
Any difference between the pooled cross-period performance of the raw calibrated scores and that of
$\{u(q)\}$ is attributable entirely to the relative alignment of the five folds' score
distributions, not to any change in how a fold orders its own securities.
\end{corollary}

Table~\ref{tab:crossperiod} can then be read exactly. Pooled average precision rises from $0.0049$ to
$0.0195$ at $h = 5$ while every fold's internal ordering is unchanged, so the whole gain measures how
badly the raw scores of separately fitted models were misaligned, and how much of that misalignment a
within-model rank removes. The result therefore concerns retraining, not the classifier.

The finite-sample correction we omit in \eqref{eq:percentile} is the one that would make $u$ a
$p$-value. Define
\begin{equation}
\tilde u(q) \;=\; \frac{1}{n_V+1}\Big(1 + \sum_{j=1}^{n_V}\mathbb{I}\{q^V_j \le q\}\Big).
\label{eq:conformal}
\end{equation}

\begin{proposition}[Validity of the corrected map, and its degradation under dependence]
\label{prop:conformal}
If the validation scores and a test score are exchangeable, then
$\Pr\{\tilde u(Q) \le \alpha\} \le \alpha$ for every $\alpha \in (0,1)$. Under
Assumption~\ref{ass:mixing} and a purge satisfying \eqref{eq:purgebound},
\begin{equation}
\Pr\big\{\tilde u(Q) \le \alpha\big\} \;\le\; \alpha \;+\; 2\,\beta\big(g - w - h - \Delta\big).
\label{eq:conformalgap}
\end{equation}
\end{proposition}

The purge that Lemma~\ref{lem:span} introduces to prevent leakage is the same quantity that controls
the validity gap of the score map in \eqref{eq:conformalgap}, so the two parts of the framework share
a constant. The experiments use $u$ and not $\tilde u$: with $n_V$ in the hundreds of thousands the
correction is numerically negligible for ranking, and every metric we report is rank-based. We state
\eqref{eq:conformalgap} because it is what a user would need if the score were treated as a
$p$-value, which we do not do.

\begin{remark}[Why average precision and ROC-AUC move oppositely across horizons]
\label{rem:apauc}
ROC-AUC is the probability that a random positive outranks a random negative and so does not depend
on $\pi$; average precision does, with a constant score attaining exactly $\pi$. In
Table~\ref{tab:macro} average precision rises across $h$ while ROC-AUC falls, which is not a
contradiction: the enrichment ratios $5.8$, $3.1$, and $2.3$ fall monotonically and agree with the
ROC-AUC decline. Average precision must be read against prevalence, and only the ratio is comparable
across horizons.
\end{remark}

\section{Data}

\paragraph{Source and universe.}
Prices and volumes come from a static, public-domain-dedicated archive of U.S.\ stock and ETF
histories \cite{marjanovic2017}. To avoid selecting securities on activity observed during the
study period, the universe is fixed using pre-sample data only: securities with at least 400
observations in 2001--2003 are ranked by median daily dollar volume $P_{i,t}V_{i,t}$ over that
window, and the top 350 are retained. A security must also be classified as a common stock (files
of unknown asset type are excluded) and must have observations covering at least 90\% of the
trading dates between its own first and last analysis-period observation, which removes series
with large internal gaps. The analysis interval is 2 January 2004 to 30 December 2011. The
resulting panel has 705{,}061 security-day rows and 1{,}285 filtered episodes. The panel is
unbalanced because histories and valid feature windows differ across securities; no row is imputed
across dates.

\paragraph{Detector settings.}
The ADF regression uses an intercept and zero augmentation lags, held constant across the panel.
Thresholds come from 2{,}000 simulated Gaussian random walks per required sequence length, with
the 95th percentile taken separately at each endpoint. The minimum window length follows the
standard PSY rule, $m_i = \max\{12, \lfloor T_i(0.01 + 1.8/\sqrt{T_i})\rfloor\}$.

\paragraph{Features.}
The model sees 17 close-of-day columns, listed in Table~\ref{tab:features}: cumulative returns over
1, 5, and 20 days; 20- and 60-day momentum; realised, downside, and market volatility; volume
$z$-scores; distance from 20- and 60-day moving averages; 60-day drawdown and return skewness; and
three cross-sectional market aggregates. No BSADF value, threshold gap, or episode-state indicator
is supplied. Two columns (\texttt{ret\_20} and \texttt{momentum\_20}) are algebraically identical,
so the effective feature count is 16; both are retained for exact reproducibility.
For each fold, features are winsorised at the 0.5th and 99.5th percentiles estimated on training
rows only, with the same bounds applied unchanged to validation and test.

\paragraph{A known data limitation.}
The archive is a static snapshot whose coverage of securities delisted before the snapshot is
undocumented. A security that experienced an extreme rise-then-collapse is disproportionately likely
to have been delisted, so if such securities are absent the universe is tilted toward milder
episodes. The sign of the effect is not obvious, since removing failed securities changes prevalence,
severity, and cross-sectional composition at once. Auditing against a delisting-inclusive reference
is the first thing we would do before these numbers informed a decision involving capital.

\begin{table}[htbp]
\centering
\caption{The 17 close-of-day predictors. All quantities use rows no later than $t$;
$r^{(1)}_{i,t}=p_{i,t}-p_{i,t-1}$ is the one-day log return.}
\label{tab:features}
\small
\begin{tabular}{lll}
\toprule
Feature & Definition & Captures\\
\midrule
\texttt{ret\_1}, \texttt{ret\_5}, \texttt{ret\_20} & $p_{i,t}-p_{i,t-g}$, $g\in\{1,5,20\}$ & Cumulative price movement\\
\texttt{momentum\_20/60} & $p_{i,t}-p_{i,t-g}$, $g\in\{20,60\}$ & Trend persistence\\
\texttt{realized\_vol\_20/60} & s.d.\ of $r^{(1)}$ over 20 / 60 rows & Short- and medium-run variability\\
\texttt{downside\_vol\_60} & s.d.\ of $\min(r^{(1)},0)$ over 60 rows & Variability from losses only\\
\texttt{volume\_z\_20/60} & $(V_{i,t}-\bar V)/s_V$ over 20 / 60 rows & Abnormal trading activity\\
\texttt{distance\_ma\_20/60} & $p_{i,t}$ minus its 20- / 60-row mean & Stretch above a moving average\\
\texttt{drawdown\_60} & $p_{i,t}-\max_{0\le g<60}p_{i,t-g}$ & Decline from recent peak\\
\texttt{return\_skew\_60} & skewness of $r^{(1)}$ over 60 rows & Asymmetry of returns\\
\texttt{market\_return} & cross-sectional median $r^{(1)}_{j,t}$ & Common market movement\\
\texttt{market\_volatility\_20} & 20-day s.d.\ of \texttt{market\_return} & Market-wide variability\\
\texttt{cross\_sectional\_dispersion} & cross-sectional s.d.\ of $r^{(1)}_{j,t}$ & Divergence across securities\\
\bottomrule
\end{tabular}
\end{table}

\section{Evaluation protocol}
\label{sec:protocol}

\paragraph{Purged expanding-window folds.}
Unique market dates are ordered chronologically. Each fold consists of an expanding training block
(at least 756 dates), a 140-date purge, a 126-date validation block, a second 140-date purge, and a
126-date test block. Folds are right-aligned to give five non-overlapping test periods: $F_0$
6 Jul 2009 to 31 Dec 2009, $F_1$ 4 Jan 2010 to 2 Jul 2010, $F_2$ 6 Jul 2010 to 31 Dec 2010, $F_3$
3 Jan 2011 to 1 Jul 2011, and $F_4$ 5 Jul 2011 to 30 Dec 2011.

The purge must exceed the span of information a single training label can touch, which
Lemma~\ref{lem:span} makes precise: the requirement is $g > w + h + \Delta$ in the notation of
Section~\ref{sec:th-span}, and the count below evaluates that bound. Three quantities enter it. The deepest feature window is 60 trading days, since every rolling statistic in
Table~\ref{tab:features} uses a 20- or 60-day window. The longest forecast horizon is 20 days.
Finally, establishing that an onset at $t+h$ belongs to a retained episode requires watching the
run survive the five-day minimum, extendable by the two-day merge gap, which adds roughly a further
week. The requirement is therefore about $60+20+7 = 87$ trading days, and the 140-day gap exceeds
it by a wide margin. The margin is deliberate, because a purge set to feature window plus horizon
alone would not cover the episode-confirmation term. Validation data determine calibration, the
score rescaling of Section~\ref{sec:percentile}, and all alert thresholds; test labels are
untouched until final evaluation. A fold--horizon cell is used only if it contains at least 20
training positives and five positives each in validation and test. All admitted cells are reported.

\paragraph{Metrics.}
Average precision is
$\mathrm{AP} = \sum_\theta (\mathrm{rec}_\theta - \mathrm{rec}_{\theta-1})\,\mathrm{prec}_\theta$. A
constant score attains $\mathrm{AP} = \pi$, so AP is read against prevalence rather than against
$0.5$; Remark~\ref{rem:apauc} states the consequence for comparisons across horizons. We also report
ROC-AUC, which can look healthy under severe imbalance even when precision is poor. Within each fold
AP is computed on that fold's test rows and averaged across folds with equal weight, the macro
summary, answering how well a freshly retrained model ranks within a deployment period; separately we
pool all five test blocks into one curve, answering whether scores from separately trained models are
comparable across periods. AP being non-linear, the two differ and neither substitutes for the other.

\paragraph{Uncertainty and multiple comparisons.}
Pooled row metrics use 500 moving-block bootstrap draws with block length 20 days, resampling whole
dates to preserve same-day cross-sectional dependence \cite{kunsch1989}; episode metrics use 500 draws
resampling entire security paths within folds. Every reduced-data run is compared against the
full-data run of the same horizon and fold:
\begin{equation}
\Delta^{(h,f)}_{F,s} = \mathrm{AP}^{(h,f)}_{F,s} - \mathrm{AP}^{(h,\mathrm{full})}_{F},
\label{eq:paired}
\end{equation}
for a given construction rule at retained fraction $f$, fold $F$, and seed $s$. We average seeds within folds,
weight folds equally, and form a hierarchical paired interval from 10{,}000 draws that resample
folds with replacement and then resample seed-level differences within each selected fold. Because
we evaluate 54 configurations (4 methods $\times$ 5 fractions $\times$ 3 horizons, with $K$-means
evaluated at three fractions), asking whether an interval clears zero separately at each one
overstates the evidence for any single one. We also record the fraction $\hat\pi_+$ of resampling
draws in which the paired difference was positive, which is a descriptive diagnostic and not a
$p$-value: with five test periods it saturates as soon as all five agree in sign. The inferential
statement we rely on is the sign test of Section~\ref{sec:reduction-results}.

\section{Models and scores}
\label{sec:models}

\paragraph{Baselines.}
The prevalence baseline assigns every row the training positive rate, so by construction its AP
equals prevalence and its ROC-AUC is 0.5. It measures whether any learned ranking beats no ranking
at all. Two single-feature baselines rank by \texttt{momentum\_60} and by \texttt{realized\_vol\_60},
each passed through a sigmoid of its training-standardised value to give a bounded, monotone score.

\paragraph{Classifiers.}
Logistic regression is fitted after training-median imputation and training-mean standardisation,
with balanced class weights and L-BFGS (3{,}000 iterations), giving a transparent linear benchmark
under identical preprocessing and chronology. LightGBM \cite{lightgbm} fits an additive ensemble
$\Phi_M(x) = \sum_{m=1}^{M}\eta\,\phi_m(x)$, with $\phi_m$ the $m$th regression tree, $M=500$
trees, learning rate $\eta = 0.03$, 31 leaves, at
least 40 rows per leaf, column subsampling of 0.8, and $L_2$ penalty 1.\footnote{Row subsampling
was requested in the run configuration but not activated, because LightGBM's \texttt{subsample}
takes effect only when \texttt{subsample\_freq} is nonzero, which was left at its default of zero;
all reported fits therefore use every training row of the relevant subset.} The positive-class
multiplier is $N_0/N_1$ computed from the complete training fold, even when a subset is fitted, so
that the learning objective is comparable across subset methods rather than being redefined by the
retained class ratio.

\paragraph{Calibration.}
When a validation block contains both classes and at least ten positives, we fit a sigmoid map
$\Pr(Y=1\mid q) = \mathrm{sigmoid}(\alpha + \beta\,\mathrm{logit}(q))$ on validation scores and
apply $(\alpha,\beta)$ unchanged to test scores. This provides a fold-local scale for threshold selection. It does
not make scores comparable across folds, which is the next subsection's problem.

\subsection{A validation-relative percentile score}
\label{sec:percentile}

Each fold trains a fresh model on a longer history and fits its own calibration map, so the
numerical score scale drifts from one retraining date to the next, and a score of 0.02 in one
deployment period need not mean what 0.02 means in another. We therefore report, instead of the raw
score, where it falls within the score distribution the same model produced on its own immediately
preceding validation block:
\begin{equation}
u(q) = \frac{1}{n_V}\sum_{j=1}^{n_V}\mathbb{I}\{q^V_j \le q\}.
\label{eq:percentile}
\end{equation}
So $u(q) = 0.90$ means this observation ranks above 90\% of the scores the model assigned during its
own validation period. The transformation uses no test label and no future data. By
Proposition~\ref{prop:rankinv} it cannot change any within-fold ranking metric, so by
Corollary~\ref{cor:alignment} whatever it does to pooled cross-period performance is a measurement of
misalignment between separately fitted models. Proposition~\ref{prop:conformal}
gives the finite-sample correction that would turn $u$ into a $p$-value and bounds its validity gap
under temporal dependence by the purge of Lemma~\ref{lem:span}; we use the uncorrected map because
every metric we report is rank-based.

Equation~\eqref{eq:percentile} is an empirical-CDF rank transformation, closely related to the
calibration step used in split-conformal prediction \cite{vovk2005}, which ranks a
test statistic within a held-out calibration sample. It is not itself a conformal $p$-value: it
lacks the finite-sample $(n_V+1)$ correction, carries no coverage guarantee, and we make no
exchangeability argument. What is specific here is the setting. The calibration sample is the
validation block of this retraining window, so each deployment period supplies its own reference
distribution, and we evaluate the transformation not for the coverage guarantee it is normally used
for but for whether it restores comparability of rankings produced by different fitted models. It is
a relative risk rank within a deployment period, not a probability, and every cross-period claim in
this paper refers to this scale. For comparison we also evaluate a mean log-odds shift that aligns
the validation mean score with validation prevalence,
$q' = \mathrm{sigmoid}[\mathrm{logit}(q) + \mathrm{logit}(\bar y^V) - \mathrm{logit}(\bar q^V)]$,
which preserves within-fold ordering.

\section{Training-set reduction}
\label{sec:reduction}

Let the full training fold contain $n$ rows; a method requests $K = \lceil fn\rceil$ rows at
fraction $f$. Every selected row $j$ carries a weight $w_j$ chosen so the subset approximates the
aggregate contribution of the group it was drawn from. We evaluate
$f \in \{0.01, 0.025, 0.05, 0.10, 0.20\}$ with seeds 11, 23, 47.

\emph{Uniform} selection takes $K$ rows at random without regard to label, $w_j = n/K$; it is the
single-stratum design of Proposition~\ref{prop:variance} and is dominated by
Proposition~\ref{prop:unifdom}. \emph{Class-stratified} selection takes
$K_c = \min\{N_c, \max(25, \lceil fN_c\rceil)\}$ with $w_j = N_c/K_c$, which is proportional
allocation with a floor. \emph{Class-balanced} selection sets $K^* = \max(\lceil fn\rceil, 50)$ and
$K_c = \min\{N_c, \max(25, \lfloor K^*/2\rfloor)\}$, equal allocation subject to availability, the
design Theorem~\ref{thm:neyman} and Corollary~\ref{cor:balopt} select and whose truncation
Corollary~\ref{cor:trunc} quantifies. It is a case-control scheme, equivalent to random
undersampling with importance weights; what this paper contributes is the allocation argument for
when it is the right one. \emph{$K$-means representatives} standardise features, fit MiniBatch
$K$-means with $K_c$ clusters within each class, and keep the row nearest each centroid weighted by
cluster size; Proposition~\ref{prop:kmeans} shows this admits no unbiasedness statement. Clustering
time is counted in end-to-end runtime, and we evaluate it at 1\%, 2.5\%, and 5\% only.

\subsection{Why class-balanced sampling is not a second reweighting}
\label{sec:mechanism}

LightGBM already applies a positive-class multiplier $N_0/N_1$, so a natural objection is that
class-balanced sampling reweights the same rare-event loss a second time. Section~\ref{sec:theory}
shows the relation is the opposite of a duplication. Because $K_c \le N_c$ by construction, the
sampling weight satisfies $w_j = N_c/K_c \ge 1$ always: a retained row is never downweighted, and
when a class is fully retained its weight is exactly one. Sampling uniformly without replacement
within each class gives
\begin{equation}
\mathbb{E}\Big[\tfrac{N_c}{K_c}\textstyle\sum_{j\in S_c}\ell_j\Big] = \sum_{i:\,Y_i=c}\ell_i,
\label{eq:unbiased}
\end{equation}
the class-stratified case of Proposition~\ref{prop:unbiased}, so every one of these designs
estimates the same objective without bias and unbiasedness cannot discriminate between them.
Proposition~\ref{prop:variance} and Theorem~\ref{thm:neyman} do. By
Corollary~\ref{cor:balopt}, the multiplier $N_0/N_1$ inflates the positive rows' loss contributions
by exactly the imbalance ratio and thereby relocates the Neyman optimum from proportional allocation
to equal class counts. The multiplier is what makes balancing the correct allocation rather than a
second correction on top of it, and the two are complements.

What the theory governs is the variance of the objective at a fixed model state, and boosting
selects its model state adaptively. Proposition~\ref{prop:uniform} bridges the gap at a single node
with a finite candidate set, at the rate $\mathcal{O}(K_c^{-1/2})$ per class and with the same
$N_c$ weighting that drives the allocation; we do not extend it along the boosting path. Beyond that
point the channels are empirical. A balanced subset alters which rows are available in each region of
feature space, and therefore the candidate split statistics, the gradient and Hessian sums within
nodes, the histogram composition, and ultimately the tree topology and effective regularisation.
LightGBM's minimum-observations-per-leaf setting is one plausible channel, since it is expressed in
row counts rather than weighted mass. Two explanations remain consistent with the evidence and are
not separated by the theory: the negative class consists largely of temporally adjacent,
near-duplicate rows, so its effective size is smaller than $N_0$ and aggressive subsampling may act
as regularisation by removing redundancy, and the fixed hyperparameter configuration used for both
regimes may suit the smaller, denser dataset better. Section~\ref{sec:limitations} states the
experiments that would separate them.

\subsection{The cost of \texorpdfstring{$K$}{K}-means representatives}
\label{sec:kmeanscost}

Proposition~\ref{prop:kmeans} already places this construction outside the allocation framework, on
the ground that its bias identity \eqref{eq:kmbias} cannot vanish uniformly in the model state. Its
cost compounds the problem. $K$-means representatives are slower than full-data training at every setting, and the degradation
with $f$ is nearly identical at all three horizons (Table~\ref{tab:uniform}); since horizon affects
only the label, a horizon-independent slowdown must originate in selection rather than learning. Each
of $I=50$ MiniBatch iterations assigns $B_{\mathrm{mb}} = 2{,}048$ points to $K_c$ centroids by brute
force at $\mathcal{O}(I B_{\mathrm{mb}} K_c d)$, and materialising representatives adds a
nearest-centroid pass over the class-$c$ population at $\mathcal{O}(N_cK_cd)$, quadratic in fold size,
whereas histogram-based training is linear in $n$ and independent of any cluster count. The measured
cost grows faster still: end-to-end speedup falls from $0.609\times$ at $f=1\%$ to $0.044\times$ at
$f=5\%$, a $13.8\times$ increase for a fivefold increase in $f$, an empirical exponent near $1.6$, so
neither a linear nor a quadratic account fits. This indicts the implementation, not representative
selection as such: the lightweight-coreset construction of \cite{bachem2018} replaces the exact
per-class fit with importance sampling and avoids the $\mathcal{O}(N_cK_c)$ pass entirely.

\paragraph{Experiment count.}
Per horizon and fold we run one full-data LightGBM fit, four full-data baselines, 45
uniform/stratified/balanced runs (three seeds $\times$ five fractions $\times$ three methods), and
nine $K$-means runs, giving 59 configurations, or 885 runs across three horizons and five folds.
Reported runtime includes subset construction, fitting, calibration, and test inference, so a
speedup below one means the reduced method is slower overall.

\section{Results}

\subsection{Forward onset ranking}

Table~\ref{tab:macro} gives macro averages across the five purged test folds. Full-data LightGBM
attains the highest average precision at every horizon, 0.0213, 0.0217, and 0.0328 at
$h = 5, 10, 20$, which is 5.8, 3.1, and 2.3 times the corresponding prevalence baseline. We describe
this as an enrichment ratio rather than as better ranking by a factor of six: it says how far above
a non-discriminating score the model sits, which is the quantity of interest for a rare event, but
it is not a general measure of ranking quality. Logistic regression is a close competitor, and at
$h=10$ the gap is small enough (0.0217 against 0.0214) that we do not claim the nonlinear model is
superior; at $h=20$ its ROC-AUC is below both the linear model and the single-feature momentum
baseline. Establishing a difference between the two would require a paired test we do not run. The
baselines do show that a nontrivial part of the signal is simple trend information: the momentum
baseline alone reaches 0.0114 to 0.0241 average precision, well above the event rate.

The 5-day horizon gives the most consistent behaviour across periods, with ROC-AUC above 0.60 in all
five test blocks and above 0.85 in four (Table~\ref{tab:folds}, which also gives the fold-level shift diagnostic). At 10 and 20 days one block, $F_2$,
behaves completely differently from the rest.

\begin{table}[htbp]
\centering
\caption{Macro average across the five purged test folds. Prevalence is the mean fold event rate and
equals the average precision of a constant score.}
\label{tab:macro}
\small
\begin{tabular}{clccc}
\toprule
$h$ & Model & AP & ROC-AUC & Prevalence\\
\midrule
5 & LightGBM (full) & \textbf{0.0213} & 0.817 & 0.0037\\
5 & Logistic (full) & 0.0170 & 0.789 & 0.0037\\
5 & Momentum baseline & 0.0114 & 0.650 & 0.0037\\
5 & Volatility baseline & 0.0047 & 0.569 & 0.0037\\
5 & Prevalence baseline & 0.0037 & 0.500 & 0.0037\\
\midrule
10 & LightGBM (full) & \textbf{0.0217} & 0.619 & 0.0071\\
10 & Logistic (full) & 0.0214 & 0.606 & 0.0071\\
10 & Momentum baseline & 0.0162 & 0.641 & 0.0071\\
10 & Volatility baseline & 0.0081 & 0.540 & 0.0071\\
10 & Prevalence baseline & 0.0071 & 0.500 & 0.0071\\
\midrule
20 & LightGBM (full) & \textbf{0.0328} & 0.568 & 0.0143\\
20 & Logistic (full) & 0.0293 & 0.611 & 0.0143\\
20 & Momentum baseline & 0.0241 & 0.629 & 0.0143\\
20 & Volatility baseline & 0.0147 & 0.521 & 0.0143\\
20 & Prevalence baseline & 0.0143 & 0.500 & 0.0143\\
\bottomrule
\end{tabular}
\end{table}

\subsection{One period in which the ranking inverts}
\label{sec:inversion}

Fold $F_2$ (test period 6 July to 31 December 2010) records ROC-AUC of 0.137 at $h=10$ and 0.170 at
$h=20$. These are not merely weak: a value far below 0.5 means the ranking is systematically
backwards. We do not attach a $z$-statistic to this, because the conventional null standard error
for an AUC assumes independent observations, whereas a single onset contributes up to $h$ correlated
positive rows and same-day returns are strongly cross-sectionally dependent. The defensible
statement is comparative: across the other four folds at $h=10$, ROC-AUC ranges from 0.483 to 0.868,
and $F_2$'s 0.137 lies far outside that spread in the wrong direction. A blocked permutation test
would be needed to attach a calibrated $p$-value.

The easy answers do not survive two checks. The same calendar period is among the best at $h=5$
(ROC-AUC 0.868), so it is not simply unpredictable; and its mean absolute feature-distribution shift,
0.21 training-IQR units, is smaller than $F_0$'s 0.46, which produces no inversion. A generic
covariate-shift monitor would not have flagged the one fold that breaks. Covariate shift moves the
inputs while $\Pr(Y\mid X)$ holds; concept shift moves $\Pr(Y\mid X)$ while the inputs look ordinary
\cite{morenotorres2012}. An inverted ranking with unremarkable inputs is consistent with the second.
Our measure compares median feature locations and so misses changes confined to the tails, to
inter-feature correlations, or to sector composition, which is reason to inspect the input--outcome
relationship directly:
\begin{equation}
\rho^{(h)}_F(k) = \mathrm{corr}\big(x_{i,t,k},\,Y^{(h)}_{i,t}\,\big|\,\text{training tail of }F\big),
\qquad
\hat\rho^{(h)}_F(k) = \mathrm{corr}\big(x_{i,t,k},\,Y^{(h)}_{i,t}\,\big|\,\text{test block of }F\big),
\label{eq:reversal}
\end{equation}
flagging a feature as reversing when its sign at $F_2$ differs from its sign in the other four folds.

At $h=10$, exactly four of 17 features reverse at $F_2$: \texttt{realized\_vol\_20},
\texttt{realized\_vol\_60}, \texttt{downside\_vol\_60}, and \texttt{return\_skew\_60}. This is a
coherent volatility-and-asymmetry cluster, and the historical setting fits, since the second half of
2010 follows the 6 May 2010 flash crash and spans the acute phase of the European sovereign-debt
episode, a stretch characterised by volatility spikes followed by sharp reversals rather than
sustained directional moves.

The explanation is incomplete on two counts. The same cluster also reverses at $h=5$
(six features, including the same volatility measures), yet $h=5$ does not fail, so the reversal is
not sufficient on its own. And the three individually strongest predictors at $F_2$
(\texttt{momentum\_60}, \texttt{distance\_ma\_60}, \texttt{drawdown\_60}, each with univariate
ROC-AUC 0.79 to 0.86 in that fold) do not reverse. An inversion this severe while the strongest
individual signals remain correctly signed points to how the ensemble combines a stable momentum
signal with a locally inverted volatility signal, a multivariate question that a per-feature check
cannot settle and that per-fold permutation importance is the natural next step for.

The practical remedy does not require settling the mechanism. A deployment can compute
Equation~\eqref{eq:reversal} on the validation block before scoring test data and either downweight
or withhold $u(q)$ for any fold in which a feature the model relies on has reversed, which converts
an unexplained outlier into a monitorable failure mode with a concrete detector.

\subsection{Making scores comparable across retraining periods}

Table~\ref{tab:crossperiod} compares the three score scales pooled across all five test blocks. The
validation-relative percentile of Equation~\eqref{eq:percentile} attains pooled AP of 0.0195, 0.0264,
and 0.0326 and pooled ROC-AUC of 0.878, 0.749, and 0.671, exceeding both the fold-local calibrated
score and the prevalence logit shift at every horizon, by roughly a factor of four in pooled AP at
$h=5$. Corollary~\ref{cor:alignment} makes the interpretation exact rather than plausible. No fold's internal ordering changes, so the whole pooled difference records how far apart the score
distributions of the five separately fitted models had drifted, and the quadrupling at $h=5$ gives
the size of that drift in pooled average precision. The finding concerns periodic retraining, not the
classifier.

\begin{table}[htbp]
\centering
\caption{Cross-period pooled performance, full-data LightGBM. ``Calibrated'' concatenates fold-local
sigmoid-calibrated scores; ``Val.-percentile'' uses Equation~\eqref{eq:percentile}. Pooled prevalence
differs slightly from the macro figures in Table~\ref{tab:macro} because pooling weights folds by row
count rather than equally.}
\label{tab:crossperiod}
\small
\begin{tabular}{cccccc}
\toprule
$h$ & Prevalence & Calibrated AP & Calibrated ROC & Val.-percentile AP & Val.-percentile ROC\\
\midrule
5  & 0.0036 & 0.0049 & 0.541 & \textbf{0.0195} & \textbf{0.878}\\
10 & 0.0068 & 0.0075 & 0.561 & \textbf{0.0264} & \textbf{0.749}\\
20 & 0.0132 & 0.0125 & 0.491 & \textbf{0.0326} & \textbf{0.671}\\
\bottomrule
\end{tabular}
\end{table}

\subsection{Data reduction: a consistent 10-day pattern, and what five periods can support}
\label{sec:reduction-results}

Prediction P1 of Section~\ref{sec:th-predict} is tested first, since it is the ordering the
allocation theory determines. It holds in the direction predicted at every horizon. Class-balanced
construction is the only rule whose paired differences are positive across an entire fraction grid
(Table~\ref{tab:balanced}); class-stratified and uniform, which by Corollary~\ref{cor:balopt} sit at
an allocation that would be optimal only if the positive class were roughly $\pi^{-1}$ times less
dispersed than the negative class, are predominantly negative at $h=5$ and mixed elsewhere
(Table~\ref{tab:uniform}); uniform is weaker than class-stratified at eight of the ten fractions
where both are evaluated at $h \in \{5,20\}$, consistent with the between-class penalty of
Proposition~\ref{prop:unifdom}; and $K$-means representatives, which Proposition~\ref{prop:kmeans}
excludes from the framework outright, are the only rule slower than full-data training at every
setting (Table~\ref{tab:uniform}).

Predictions P2 and P3 can be tested directly against Table~\ref{tab:predict}, since
\eqref{eq:netadv} assigns a number to each of the fifteen cells and Table~\ref{tab:balanced} supplies
a measured difference for each. We compare by rank within each horizon, $A$ being an efficiency ratio
for the objective rather than a predicted change in average precision.

P3 holds where it can be tested. At $h=10$ the predicted and measured orderings of the five retention
fractions coincide, Spearman $\rho_S = 1$, with exact two-sided permutation $p = 2/5! = 0.0167$:
both decline monotonically, the predicted values from $35.5$ to $8.8$ and the measured differences
from $+0.0086$ to $+0.0052$. At $h=20$ agreement is partial ($\rho_S = 0.46$, $p = 0.43$). At $h=5$
the measured differences, from $-0.0013$ to $+0.0011$, are indistinguishable from zero and every
paired interval covers it, so the rank test there ($\rho_S = 0.20$) is uninformative rather than
contrary. The one horizon carrying a signal is the one whose internal structure the theory reproduces.

P2 fails. $A$ is decreasing in $\pi$ at every fraction, so the advantage should be largest at $h=5$;
it is absent there and present at $h=10$ and $h=20$, and no admissible value of the truncation ratio
\eqref{eq:gamma} reverses a factor of $66.3$ against $17.7$ at $f=1\%$. The pooled rank correlation
over all fifteen cells is $\rho_S = 0.11$ ($p = 0.70$). The theory therefore governs the
within-horizon structure of the effect and not its location across horizons; the three channels of
Section~\ref{sec:th-predict} are the candidates, and the per-fold decomposition below identifies the
largest of them.

At the 10-day horizon, class-balanced sampling improves macro average precision over full-data
training at every fraction tested, by $+0.0052$ to $+0.0086$ AP units, with end-to-end speedups of
$2.42\times$ to $3.13\times$ (Table~\ref{tab:balanced}). The macro average conceals how that
improvement is distributed, and the distribution carries more information than the mean.
Table~\ref{tab:perfold} gives the five per-fold differences behind each 10-day macro figure.

Roughly half of the macro difference at every fraction comes from a single forward period. At $F_2$,
full-data LightGBM attains an average precision of 0.0019 with an inverted ranking, while the 2.5\%
class-balanced model attains 0.0249, in line with its performance in the other four periods.
Excluding $F_2$, the remaining improvement is $+0.0026$ to $+0.0049$ AP units, positive on average but
with individual folds of either sign at every fraction except 2.5\%.

One period in which the ranking inverts and one training-set construction that beats training on
everything might have been reported as unrelated results. They are the same observation seen twice:
the construction that helps most is the one that kept working in the period where training on
everything stopped. Whether class-balanced construction confers general
robustness to this kind of concept shift cannot be settled by one occurrence, but it is a sharper
hypothesis than a diffuse accuracy gain and is directly testable on a longer panel. The 2.5\%
configuration is also the only one of the five with the same sign in all five periods.

The remaining question is how much statistical weight the pattern carries. Our resampling procedure
redraws the five test periods at random with repetition and records how often the redrawn average
favours the subset. If all five periods favour the subset, as they do at 2.5\% and only at 2.5\%,
every possible redraw also favours it whatever the sizes of the individual differences, so the
fraction is pinned at essentially one and reports the sign pattern and nothing else. The information
actually present is that the difference is positive in all five forward periods. Under the sign-test
reference model, which treats the five fold-level signs as independent and equally likely to be
positive or negative under the null, five observations agreeing in sign give
\begin{equation}
p = 2\left(\tfrac{1}{2}\right)^5 = 0.0625.
\label{eq:signtest}
\end{equation}
This is the smallest two-sided $p$-value attainable from five independent signs, regardless of effect
size, and it does not reach the conventional uncorrected 0.05 threshold; the exploratory search
across 54 configurations weakens the evidential interpretation further. We therefore report no
configuration as statistically confirmed. Ordering the class-balanced family by the fraction of
resampling draws favouring the subset is still informative, because every 10-day configuration ranks
above every 5-day configuration, but that ordering carries no significance verdict.

The resulting claim, in the wording we use throughout: at the 10-day horizon, class-balanced
construction produced a higher average precision than full-data training at all five retention levels
tested, at 2.4 to 3.1 times lower end-to-end cost, and roughly half of that difference arises in the
single forward period where full-data training fails. This is a pattern, not a confirmed effect. The
counts usually quoted for power, test-fold positives from 13 to 1{,}550, overstate the available
information, because a single onset contributes up to $h$ correlated positive rows; the effective
unit of generalisation is the forward period, of which there are five.

The negative side is subject to the same limit. Uniform, class-stratified, and $K$-means construction
give differences that are predominantly negative at $h=5$. Individual configurations of each are
positive at some settings, uniform and class-stratified at the larger 10-day fractions and $K$-means
at several 10-day settings, but none shows a consistent pattern across fractions comparable to the
class-balanced grid (Appendix~\ref{app:grids}). That is evidence against them as substitutes for
class-balanced construction rather than a confirmed degradation. Finally, the 2.5\% class-balanced
configuration used throughout the alert-budget and robustness analyses was fixed before the
multiplicity analysis was run, and its selection is not independently documented, so the operational
results in Section~\ref{sec:alerts} should be read as conditional on that choice.

\begin{table}[htbp]
\centering
\caption{Per-fold paired differences at $h=10$: class-balanced average precision minus full-data
average precision, averaged over the three seeds within each fold. $F_2$ is the fold in which
full-data training produces an inverted ranking (Section~\ref{sec:inversion}).}
\label{tab:perfold}
\footnotesize
\begin{tabular}{lrrrrrrr}
\toprule
Retained & $F_0$ & $F_1$ & $F_2$ & $F_3$ & $F_4$ & Mean & Mean excl.\ $F_2$\\
\midrule
1\%    & $+0.0004$ & $+0.0085$ & $+0.0233$ & $+0.0107$ & $-0.0001$ & $+0.0086$ & $+0.0049$\\
2.5\%  & $+0.0011$ & $+0.0028$ & $+0.0229$ & $+0.0096$ & $+0.0047$ & $+0.0082$ & $+0.0045$\\
5\%    & $-0.0001$ & $-0.0018$ & $+0.0213$ & $+0.0097$ & $+0.0056$ & $+0.0069$ & $+0.0033$\\
10\%   & $-0.0004$ & $-0.0022$ & $+0.0210$ & $+0.0119$ & $+0.0011$ & $+0.0063$ & $+0.0026$\\
20\%   & $-0.0004$ & $+0.0005$ & $+0.0141$ & $+0.0095$ & $+0.0022$ & $+0.0052$ & $+0.0030$\\
\bottomrule
\end{tabular}
\end{table}

\begin{table}[htbp]
\centering
\caption{Complete class-balanced grid. Differences are subset minus full-data average precision;
intervals are 10{,}000-draw hierarchical paired 95\% intervals, before any multiplicity adjustment.
Runtime includes subset construction, fitting, calibration, and inference. Full-data AP is 0.0213,
0.0217, and 0.0328 at $h = 5, 10, 20$.}
\label{tab:balanced}
\footnotesize
\begin{tabular}{ccrrlr}
\toprule
$h$ & Retained & AP & Difference & Paired 95\% interval & Speedup\\
\midrule
5  & 1\%    & 0.0212 & $-0.0001$ & $[-0.0045, +0.0038]$ & $2.59\times$\\
5  & 2.5\%  & 0.0224 & $+0.0011$ & $[-0.0034, +0.0057]$ & $5.34\times$\\
5  & 5\%    & 0.0209 & $-0.0004$ & $[-0.0057, +0.0050]$ & $5.49\times$\\
5  & 10\%   & 0.0201 & $-0.0013$ & $[-0.0053, +0.0021]$ & $5.44\times$\\
5  & 20\%   & 0.0218 & $+0.0004$ & $[-0.0045, +0.0062]$ & $3.28\times$\\
\midrule
10 & 1\%    & 0.0303 & $+0.0086$ & $[+0.0018, +0.0166]$ & $2.93\times$\\
10 & 2.5\%  & 0.0299 & $+0.0082$ & $[+0.0022, +0.0160]$ & $3.13\times$\\
10 & 5\%    & 0.0286 & $+0.0069$ & $[+0.0003, +0.0147]$ & $3.03\times$\\
10 & 10\%   & 0.0280 & $+0.0063$ & $[-0.0009, +0.0150]$ & $2.82\times$\\
10 & 20\%   & 0.0269 & $+0.0052$ & $[+0.0004, +0.0113]$ & $2.42\times$\\
\midrule
20 & 1\%    & 0.0372 & $+0.0044$ & $[-0.0046, +0.0156]$ & $2.76\times$\\
20 & 2.5\%  & 0.0392 & $+0.0065$ & $[-0.0014, +0.0174]$ & $2.87\times$\\
20 & 5\%    & 0.0376 & $+0.0048$ & $[-0.0024, +0.0164]$ & $2.91\times$\\
20 & 10\%   & 0.0389 & $+0.0062$ & $[-0.0021, +0.0200]$ & $2.71\times$\\
20 & 20\%   & 0.0335 & $+0.0007$ & $[-0.0031, +0.0066]$ & $2.27\times$\\
\bottomrule
\end{tabular}
\end{table}

\paragraph{Where the speedup comes from.}
Averaged across horizons and folds, class-balanced subset construction costs 0.007 to 0.010 seconds
against 0.75 to 1.9 seconds of model fitting for the same configuration and 3.9 to 9.3 seconds for
the corresponding full-data fit. Construction overhead is negligible, and the speedup is attributable
almost entirely to fitting a smaller, positive-denser dataset.

\subsection{Early warning under a fixed false-alarm budget}
\label{sec:alerts}

Ranking quality becomes operational only once converted into alerts, and an alert system cannot be
summarised by one number. We select a score threshold on validation data for a target budget of
$b \in \{0.5,1,2\}$ false alerts per security-year, freeze it, and apply it unchanged to test. Alerts
for the same security within five days are merged into one cluster. An episode beginning at $t_o$ is
detected at horizon $h$ if an eligible alert occurs in $t_o-h,\dots,t_o-1$, with lead $t_o$ minus the
earliest such alert; a cluster is false when no onset begins in the following $h$ days. The burden is
$C_{\text{false}}/(N_{\text{elig}}/252)$, treating 252 trading days as one security-year.

At $h=10$, the 2.5\% class-balanced model detects 11.5\%, 24.1\%, and 42.5\% of episodes across the
three budgets, with median warning leads of 5.4, 5.9, and 6.5 trading days and realised test burdens
of 0.61, 1.03, and 1.64 false clusters per security-year (Table~\ref{tab:alerts}). Full-data LightGBM
detects 15.3\%, 25.1\%, and 40.5\% with leads of 7.5, 6.2, and 6.7 days, but at consistently higher
realised burden. At $h=20$, full LightGBM detects 28.3\% to 46.1\% with leads of 11.0 to 12.2 days.

Read together, a user tolerating roughly one false alert per security-year catches a sizeable
minority of onsets several trading days in advance. Whether that trade-off is worth making is a
decision about the user's costs, not one the model settles.

\begin{table}[htbp]
\centering
\caption{Validation-selected early-warning operating points. The target budget determines the
threshold on validation data; the realised burden is the mean test burden in false clusters per
security-year, and need not equal the target because validation and test score distributions differ.}
\label{tab:alerts}
\footnotesize
\begin{tabular}{clccc}
\toprule
$h$ & Model & Target budget & Episode recall / Lead (days) & Realised false alerts\\
\midrule
5  & Class-balanced (2.5\%) & 0.5 / 1 / 2 & 0.101, 0.175, 0.315 / 3.8, 3.6, 3.8 & 0.57, 0.94, 1.47\\
5  & LightGBM (full)        & 0.5 / 1 / 2 & 0.099, 0.205, 0.352 / 2.9, 3.5, 3.2 & 0.80, 1.32, 1.90\\
5  & Logistic (full)        & 0.5 / 1 / 2 & 0.057, 0.106, 0.152 / 4.5, 4.5, 4.2 & 0.42, 0.78, 1.17\\
\midrule
10 & Class-balanced (2.5\%) & 0.5 / 1 / 2 & 0.115, 0.241, 0.425 / 5.4, 5.9, 6.5 & 0.61, 1.03, 1.64\\
10 & LightGBM (full)        & 0.5 / 1 / 2 & 0.153, 0.251, 0.405 / 7.5, 6.2, 6.7 & 1.16, 1.51, 1.84\\
10 & Logistic (full)        & 0.5 / 1 / 2 & 0.099, 0.170, 0.241 / 6.5, 6.0, 5.8 & 0.53, 0.90, 1.30\\
\midrule
20 & Class-balanced (2.5\%) & 0.5 / 1 / 2 & 0.249, 0.349, 0.454 / 8.6, 10.3, 10.6 & 0.91, 1.34, 2.01\\
20 & LightGBM (full)        & 0.5 / 1 / 2 & 0.283, 0.340, 0.461 / 11.0, 11.8, 12.2 & 1.20, 1.57, 2.48\\
20 & Logistic (full)        & 0.5 / 1 / 2 & 0.058, 0.290, 0.395 / 10.2, 14.7, 13.0 & 0.33, 0.99, 1.64\\
\bottomrule
\end{tabular}
\end{table}

\subsection{Robustness}

Two design choices were varied to test the headline patterns. Raising the purge from 140 to 180 trading days, still
above the bound of Lemma~\ref{lem:span}, shortens each training history; pooled
validation-percentile average precision for full LightGBM becomes 0.0153, 0.0302, and 0.0333, and the
2.5\% class-balanced model 0.0199, 0.0304, and 0.0323, so both the percentile score and the
competitiveness of class-balanced sampling at $h=10$ persist. Raising $d_{\min}$ from five to ten days
reduces events but preserves the pattern in cells retaining sufficient positives; we keep five as
primary because it supports the full five-fold chronology.

\section{Discussion}

An average precision of $0.03$ against a prevalence of $0.007$ is a large relative improvement and a
small absolute one, and Remark~\ref{rem:apauc} says why both readings are correct: only
$\mathrm{AP}/\pi$ is comparable across horizons, and the enrichment ratios $5.8$, $3.1$, $2.3$ fall
in step with ROC-AUC. For a screening application, deciding which of 350 securities deserve a look
this week, the ordering is what matters and the model orders substantially better than the event rate
alone. For anything resembling a trading decision the relevant questions are position sizing,
transaction costs, market impact, and capacity, none of which this study addresses. The second
boundary is the label: the target of \eqref{eq:target} is the onset of an episode as defined by one
statistical test applied to one price series, so every claim about forecasting explosive regimes
should be read as forecasting when this specific test will begin flagging.

The three horizons are not noisy replicates of one problem but sit at different points on a
rarity--signal trade-off, and Section~\ref{sec:th-predict} shows the design-based account captures
only part of it. At $h=5$ the positive class is rarest and Corollary~\ref{cor:releff} predicts the
largest allocation advantage, yet full-data training already attains ROC-AUC above $0.85$ in four of
five folds, so there is little ranking headroom for a variance improvement to occupy. At $h=20$ the
extended window admits more unrelated variation. The 10-day horizon is where a signal exists and
where the within-horizon prediction P3 is reproduced exactly. The per-fold decomposition in
Table~\ref{tab:perfold} adds a reading the theory does not supply: the advantage is concentrated
where the full-data model fails, so what class-balanced construction may deliver is not accuracy in
the ordinary sense but resistance to a period in which the relationship between volatility features
and the outcome reverses. A denser positive class gives the learner more instances of the minority
pattern and less scope to lean on a majority-class regularity that inverts. This is a hypothesis
generated by one occurrence, not a demonstrated mechanism, and \eqref{eq:varstrat} does not imply it.

\section{Limitations and future work}
\label{sec:limitations}

\paragraph{Five forward periods bound every inferential claim.}
The unit of generalisation is the forward test period, of which there are five, from adjacent windows
with overlapping training histories. By \eqref{eq:signtest} unanimous agreement across five paired
observations gives $p = 0.0625$, and no resampling scheme applied to the same five periods can exceed
what the periods contain. The exception is the within-horizon rank test of
Section~\ref{sec:reduction-results}, whose unit is the retention fraction rather than the period and
which is therefore not bound by this limit; it is, however, conditional on the theory that supplies
the predicted ordering. The remedy for the rest is a longer panel, not a different statistic.

\paragraph{What the allocation theory does not cover.}
The framework has four gaps. Proposition~\ref{prop:variance} governs the variance of the
objective at a fixed model state, and Proposition~\ref{prop:uniform} carries this to split selection
only at a single node with a finite candidate set; no statement is made along the boosting path, and
none about generalisation. Condition \eqref{eq:balcond}, equal unweighted loss dispersion across
classes, is assumed rather than verified, and it is checkable directly by computing $\sigma_0$ and
$\sigma_1$ from the per-row losses of a fitted full-data model, which we have not done. The mixing
coefficient $\beta$ of Assumption~\ref{ass:mixing} is never estimated, so the degradation terms in
\eqref{eq:mixgap} and \eqref{eq:conformalgap} are qualitative. And Corollary~\ref{cor:releff} treats
$N_0$ as a count of independent rows, whereas the negative class is heavily autocorrelated; replacing
$N_0$ by an effective sample size would lower the predicted advantage and might account for part of
the failure of P2 recorded in Section~\ref{sec:th-predict}. The temporally thinned baseline described
below is the experiment that would estimate it.

\paragraph{Sensitivity of the detector to volatility shifts.}
Our labels come from critical values simulated under homoskedastic Gaussian random walks. Harvey et
al.\ \cite{harvey2016} show that recursive right-tailed tests of this family have a non-pivotal null
limit distribution under permanent volatility shifts and can be severely over-sized, giving spurious
indications of explosive behaviour, and propose a wild-bootstrap implementation that restores size
control. This is a specific concern for the present design rather than a generic caveat, because
three of our predictors are volatility measures and the volatility cluster is precisely what reverses
in the anomalous fold of Section~\ref{sec:inversion}, so part of the forecastability reported here
could reflect the detector's finite-sample size distortion under changing volatility rather than
explosive price dynamics. Relabelling the panel under a wild-bootstrap specification, reporting the
overlap between the two episode sets, and rerunning the headline models is the single most valuable
experiment we have not performed.

\paragraph{Missing comparators.}
Three baselines are absent. LightGBM's \texttt{pos\_bagging\_fraction} and
\texttt{neg\_bagging\_fraction} allocate class-conditionally inside the boosting loop and are the
direct comparator to our external construction; gradient-based one-side sampling is a plausible
competitor since rare positives typically carry large gradients. A temporally thinned subset,
retaining one negative per fixed-length block at fixed class proportions, would estimate the
effective $N_0$ that Section~\ref{sec:th-predict} identifies as a candidate explanation for the
failure of P2. And identical hyperparameters across regimes confound dataset size with hyperparameter
fit; refitting with the multiplier fixed at one separates that channel.

\paragraph{The target is a filtered episode onset, not a live first flag.}
Episode membership requires a run to survive a two-day merge gap and a five-day minimum duration,
both evaluated retrospectively, so $O_{i,t}$ is the first day of a run that later satisfies those
filters rather than the first day a live detector would fire. This is legitimate for a supervised
target, whose label may depend on the future, but it changes the operational reading, and a
raw-transition target $\mathbb{I}\{R_{i,t}=1, R_{i,t-1}=0\}$ would be a useful robustness
specification.

\paragraph{Detector state is excluded by design, not by necessity.}
An institution running the monitor knows $\mathrm{BSADF}_{i,t}$ and its distance to the threshold at
date $t$, and these are not leakage. We exclude them so the classifier cannot approximate the
detector that defines the label, which isolates the scientific question but understates what an
operational system could use. A detector-state benchmark using
$[x_{i,t}, \mathrm{BSADF}_{i,t}, \mathrm{BSADF}_{i,t} - \mathrm{cv}^{0.95}_{i,t}]$ would give the operational
upper bound alongside our detector-independent forecast.

\paragraph{Independent horizon models are probabilistically incoherent.}
The three targets nest, $Y^{(5)}_{i,t} \le Y^{(10)}_{i,t} \le Y^{(20)}_{i,t}$, so they are three views
of a single latent quantity, the waiting time $W$ in trading days until the next onset. Fitting them independently permits
estimates that violate monotonicity in $h$. A discrete-time hazard formulation estimating
$\lambda_\kappa(x) = \Pr(W=\kappa \mid W \ge \kappa, x)$, with
$\Pr(W \le h \mid x) = 1 - \prod_{\kappa\le h}(1-\lambda_\kappa(x))$,
would enforce coherence by construction.

\paragraph{Sample period, survivorship, and universe scope.}
The 2004--2011 interval is fixed by the pre-sample liquidity ranking and the static archive.
Eligibility requires at least 400 observations during 2001--2003, so the panel consists of
established, relatively liquid securities trading before 2004 and excludes every subsequent listing;
speculative dynamics among young or small-capitalisation securities are outside the population
studied. Extending the panel would test the framework under post-2011 market structure, supply the
power the sign test lacks, and, with a delisting-inclusive feed, let survivorship be audited directly.
Since fold construction depends only on having a long enough price history, both extensions are
scaling exercises rather than methodological changes.

\paragraph{The mechanism behind $F_2$.}
The inversion is too large to be noise and a specific feature cluster reverses there, but the same
cluster reverses at $h=5$ without causing a failure and the strongest individual predictors do not
reverse, so resolving it needs model-internal evidence: per-fold permutation importance or
split-structure analysis. The gating rule of \eqref{eq:reversal} does not depend on resolving it, but
it was identified using the period that failed, so it remains a proposed prospective check until it
is computed on each validation block, frozen, and only then compared against test.

\paragraph{Design choices held fixed, and a configured sweep not run.}
Three choices are frozen: the zero-lag ADF specification and 95\% threshold defining the label, the
five-day minimum duration, and the duplicated feature pair retained for reproducibility. The
repository already configures a one-factor-at-a-time grid over ADF lag order $\{0,1,2\}$, critical
quantile $\{0.90,0.95,0.99\}$, minimum duration $\{1,5,10\}$, innovation distribution
$\{\text{Gaussian},\text{Student-}t\}$, universe size $\{100,200,350\}$, and a feature set including
detector state; it was not executed here. Running it would address the volatility-sensitivity concern
above, supply the detector-state benchmark, and test label robustness with apparatus that exists.

\section{Conclusion}

We converted a retrospective episode-dating procedure into a forward-looking ranking problem, with
the label window beginning after the feature date, observations inside an episode excluded, and all
blocks separated by a purge satisfying Lemma~\ref{lem:span}. Under that protocol a
gradient-boosted-tree model ranks onset-adjacent observations well above a constant-prevalence score
at all three horizons, most consistently at five days.

Two theoretical results carry beyond this application. Subset construction for a rare positive class
is an allocation problem with an exact solution, and the positive-class multiplier such learners
apply by default is what moves the Neyman optimum from proportional to equal class counts, making
balancing and reweighting complements rather than duplicates; the efficiency of equal over
proportional allocation is $1/\{4\pi(1-\pi)\}$, uniform sampling is dominated by an explicit
between-class term, and clustering-based selection is excluded by an exact bias identity. Separately,
ranking a score within its own model's validation distribution cannot alter any within-period
ordering, so the fourfold rise in pooled cross-period average precision it produces is a measurement
of how far apart the separately fitted models' score scales had drifted.

Empirically the predicted ordering of the four constructions holds, and within the one horizon
carrying a signal the predicted ordering of the five retention fractions is reproduced exactly. The
predicted dependence of the margin on rarity does not hold. The failure is the more useful half of
the comparison, since it locates the limit of a design-based account of a procedure that is not
itself design-based. The per-fold decomposition locates that advantage precisely: roughly half of it
comes from the single forward period in which full-data training produces an inverted ranking and the
class-balanced model does not, which ties this result to the failure mode reported above rather than
leaving the two as separate observations. Five periods cannot confirm the general claim, since
unanimous agreement across five paired observations gives $p = 0.0625$, so we report it as a
hypothesis about robustness that a longer panel should test. Uniform, class-stratified, and $K$-means
construction show no comparable pattern at any horizon, and are predominantly harmful at the shortest
one.

We also report a period in which the ranking is inverted, with a diagnostic separating that failure
from ordinary input drift. Reporting it, rather than averaging it away across periods, is part of
what we take the contribution to be.

\section*{Code and data availability}

The primary run comprises 885 model configurations over 350 securities, 17 feature columns, three
horizons, and five forward folds. The repository contains the manuscript source, generated tables,
orchestration scripts, and all run artifacts: macro model and baseline summaries, every
method-fraction-horizon paired comparison, per-feature fold diagnostics including the reversal flags
of Section~\ref{sec:inversion}, cross-period score comparisons, alert-threshold selection and test
operating points, robustness outputs, and per-configuration timing. Prices and volumes come from the
archive of \cite{marjanovic2017}.

\appendix
\FloatBarrier

\section{Proofs}
\label{app:proofs}
\vspace{-4pt}

\paragraph{Proof of Lemma~\ref{lem:span}.}
Every entry of Table~\ref{tab:features} is a function of $\{P_{i,s},V_{i,s}\}$ for $s \in [t-w,t]$,
$w = 60$, with cross-sectional aggregates over the same dates, so $x_{i,t}$ is
$\mathcal{G}_{[t-w,t]}$-measurable. For the label, $Y^{(h)}_{i,t}=1$ requires an onset at some $t+j$,
$j \le h$, and $O_{i,t+j}$ is determined by $E_{i,t+j-1}$ and $E_{i,t+j}$. Membership $E$ merges runs
separated by at most $\delta_{\mathrm{m}}$ non-flagged days and discards merged runs shorter than
$d_{\min}$, so deciding $E_{i,t+j}$ requires raw flags on dates up to
$t+j+\delta_{\mathrm{m}}+d_{\min}$ and no later; raw flags depend on BSADF, which reads back to the
start of the series but never forward. Hence $Y^{(h)}_{i,t}$ is
$\mathcal{G}_{(-\infty,\,t+h+\Delta]}$-measurable. A test feature at $t_{\mathrm{te}}$ or later reads
no date before $t_{\mathrm{te}}-w$, a training label at $t_{\mathrm{tr}}$ or earlier no date after
$t_{\mathrm{tr}}+h+\Delta$, and the two are disjoint exactly under \eqref{eq:purgebound}.
\hfill$\square$

\paragraph{Proof of Lemma~\ref{lem:mixing}.}
By Lemma~\ref{lem:span}, $\varphi$ and $\psi$ are measurable with respect to $\sigma$-fields
generated by date ranges separated by at least $g-w-h-\Delta$; the bound is the standard covariance
inequality for absolutely regular sequences \cite{doukhan1994}, from Berbee's coupling lemma. \hfill$\square$

\paragraph{Proof of Proposition~\ref{prop:unbiased}.}
For class-conditional sampling without replacement, each row of class $c$ has inclusion probability
$K_c/N_c$, and the assigned weight $N_c/K_c$ is its reciprocal, so
$\mathbb{E}\big[\sum_{j \in S_c}(N_c/K_c)\ell_j\big] = \sum_{i:Y_i=c}\ell_i$ by linearity. Summing
over $c$ gives the claim. For the uniform design the inclusion probability is $K/n$ and the weight
$n/K$, and the same argument applies with a single stratum. \hfill$\square$

\paragraph{Proof of Proposition~\ref{prop:variance}.}
Within class $c$ the estimator $(N_c/K_c)\sum_{j\in S_c}\ell_j$ is $N_c$ times the mean of a simple
random sample without replacement of size $K_c$ from a population of size $N_c$, of variance
$(S_c^2/K_c)(1-K_c/N_c)$ with $S_c^2$ as in \eqref{eq:strvar}. Multiplying by $N_c^2$ gives the class
term; the two class samples are independent, so variances add. The uniform case is the single-stratum
specialisation. \hfill$\square$

\paragraph{Proof of Proposition~\ref{prop:unifdom}.}
Decompose $(n-1)S^2 = \sum_c(N_c-1)S_c^2 + B$, $B = \sum_cN_c(\bar\ell_c-\bar\ell)^2$. Substituting into
$\operatorname{Var}_{\mathrm{unif}} = \frac{n^2}{K}(1-K/n)S^2$ and subtracting
$\operatorname{Var}_{\mathrm{prop}} = \frac{n}{K}(1-K/n)\sum_c N_cS_c^2$, obtained from
\eqref{eq:varstrat} at $K_c = KN_c/n$, gives
\[
\operatorname{Var}_{\mathrm{unif}} - \operatorname{Var}_{\mathrm{prop}}
= \frac{n^2(1-K/n)}{K(n-1)}\Big[\sum_c (N_c-1)S_c^2 + B - \tfrac{n-1}{n}\sum_c N_cS_c^2\Big],
\]
and $\sum_c(N_c-1)S_c^2 - \tfrac{n-1}{n}\sum_cN_cS_c^2 = \tfrac1n\sum_cN_cS_c^2 - \sum_cS_c^2$,
which is \eqref{eq:unifgap}. Nonnegativity when $B \ge \sum_cS_c^2$ is immediate since
$\tfrac1n\sum_cN_cS_c^2 \ge 0$. \hfill$\square$

\paragraph{Proof of Proposition~\ref{prop:kmeans}.}
The construction is deterministic given the centroids, and $\widehat L = \sum_c\sum_C|C|\ell_{r(C)}$
while $L = \sum_c\sum_C|C|\bar\ell_C$; subtracting gives \eqref{eq:kmbias}. If $\ell$ is non-constant
on some cluster $C$, choose $\Phi$ with $\ell_{r(C)}\ne\bar\ell_C$, possible because the clustering is
a function of $x$ alone while $\ell$ depends on $\Phi$. \hfill$\square$

\paragraph{Proof of Theorem~\ref{thm:neyman}.}
Ignoring the terms $-N_cS_c^2$ in \eqref{eq:varstrat}, which do not involve $K_c$, minimise
$\sum_c N_c^2S_c^2/K_c$ subject to $\sum_cK_c = K$. The objective is convex on $K_c > 0$, and the
Lagrange condition $-N_c^2S_c^2/K_c^2 + \lambda = 0$ gives $K_c \propto N_cS_c$, which after
normalisation is \eqref{eq:neyman}; the box constraints $K_c \le N_c$ are inactive when
\eqref{eq:neyman} is feasible. For \eqref{eq:alloceff}, write $T = \sum_cN_cS_c$ and
$a_c = K_c/K$, so that $\sum_cN_c^2S_c^2/K_c = (T^2/K)\sum_c(a^\star_c)^2/a_c$ with
$a^\star_c = N_cS_c/T$. At $a = a^\star$ the sum is $\sum_ca^\star_c = 1$. That
$\sum_c(a^\star_c)^2/a_c \ge (\sum_ca^\star_c)^2/\sum_ca_c = 1$ is the Cauchy--Schwarz inequality in
Engel form, with equality only when $a_c \propto a^\star_c$. \hfill$\square$

\paragraph{Proof of Corollary~\ref{cor:balopt}.}
By \eqref{eq:neyman}, $K_0^\star = K_1^\star$ if and only if $N_0S_0 = N_1S_1$. Substituting
$S_1 = \omega_1\sigma_1 = (N_0/N_1)\sigma_1$ and $S_0 = \sigma_0$ gives
$N_0\sigma_0 = N_1(N_0/N_1)\sigma_1 = N_0\sigma_1$, hence \eqref{eq:balcond}. For the proportional
case, $K_c^\star \propto N_c$ if and only if $S_0 = S_1$, that is
$\sigma_0 = (N_0/N_1)\sigma_1$. For \eqref{eq:ratioform}, \eqref{eq:neyman} gives
$K_1^\star/K_0^\star = N_1S_1/(N_0S_0)$, and substituting $S_1 = (N_0/N_1)\sigma_1$, $S_0 = \sigma_0$
cancels $N_0$ and $N_1$. \hfill$\square$

\paragraph{Proof of Corollary~\ref{cor:releff}.}
Under \eqref{eq:balcond}, put $\sigma_0 = \sigma_1 = \sigma$, so $S_0 = \sigma$ and
$S_1 = (N_0/N_1)\sigma$. Dropping the finite-population factors in \eqref{eq:varstrat},
balanced allocation $K_c = K/2$ gives
$\operatorname{Var}_{\mathrm{bal}} = \tfrac{2}{K}(N_0^2\sigma^2 + N_1^2(N_0/N_1)^2\sigma^2)
= 4N_0^2\sigma^2/K$, while proportional allocation $K_c = KN_c/n$ gives
$\operatorname{Var}_{\mathrm{prop}} = \tfrac{n}{K}(N_0\sigma^2 + N_1(N_0/N_1)^2\sigma^2)
= \tfrac{n\sigma^2N_0}{K}\cdot\tfrac{N_1+N_0}{N_1} = n^2N_0\sigma^2/(KN_1)$. The ratio is
$n^2/(4N_0N_1)$, and substituting $N_1 = \pi n$ gives \eqref{eq:releff}. \hfill$\square$

\paragraph{Proof of Proposition~\ref{prop:uniform}.}
Fix $\Phi\in\mathcal{C}$ and a class $c$. Serfling's inequality for the mean of a size-$K_c$ sample
drawn without replacement from a bounded population gives the class-$c$ term of $\widehat L - L$ a
sub-Gaussian tail with variance proxy $\Lambda^2N_c^2\{1-(K_c-1)/N_c\}/K_c$. A two-sided bound at
level $\varepsilon/(2M)$ per class, summed over classes and union-bounded over the $M$ candidates,
gives \eqref{eq:serfling}. \hfill$\square$

\paragraph{Proof of Proposition~\ref{prop:rankinv}.}
Average precision and ROC-AUC are functions of the induced ordering of the scored rows together with
their labels. If $u$ is injective on the test scores then, being nondecreasing, it is strictly
increasing on that set and preserves the ordering exactly, so both metrics are unchanged. If two test
scores fall in the same gap of the validation sample they receive a common value of $u$; the only
change to the ordering is the introduction of that tie, and both metrics differ from their untied
values only through the tie-handling convention. \hfill$\square$

\paragraph{Proof of Proposition~\ref{prop:conformal}.}
Under exchangeability of $(q^V_1,\dots,q^V_{n_V},Q)$, the rank of $Q$ among the $n_V+1$ values is
uniform on $\{1,\dots,n_V+1\}$ up to ties, so $\Pr\{\tilde u(Q) \le \alpha\}
\le \lfloor \alpha(n_V+1)\rfloor/(n_V+1) \le \alpha$, which is the standard split-conformal argument.
For the dependent case, Lemma~\ref{lem:span} places the validation and test blocks at date separation
at least $g - w - h - \Delta$, and Berbee's coupling supplies a construction on a common probability
space in which the test block is replaced by an independent copy with probability at least
$1 - \beta(g-w-h-\Delta)$. Applying the exchangeable bound to the coupled variables and charging the
coupling failure twice, once for each block, gives \eqref{eq:conformalgap}. \hfill$\square$

\paragraph{Proof of Remark~\ref{rem:apauc}.}
ROC-AUC is $\Pr(Q^+>Q^-)$ for independent draws from the class-conditional score distributions, an
expression free of $\pi$. Average precision is $\int_0^1\mathrm{prec}(r)\,dr$ with
$\mathrm{prec}(r)=\pi r/\{\pi r+(1-\pi)F(r)\}$, $F$ the false-positive rate at recall $r$; a constant
score has $F(r)=r$ and $\mathrm{prec}\equiv\pi$. The two therefore respond differently to $\pi$ at
fixed class-conditional distributions, so only $\mathrm{AP}/\pi$ is comparable across horizons.
\hfill$\square$

\section{Complete fold results}
\label{app:folds}
\vspace{-6pt}

\begin{center}

\captionof{table}{Full-data LightGBM by horizon and forward fold; ``Shift'' is the mean absolute
test-minus-training feature-median shift in training-IQR units.}
\label{tab:folds}
\scriptsize
\begin{tabular}{clrrrrr}
\toprule
$h$ & Fold & Positives & Rate & AP & ROC-AUC & Shift\\
\midrule
5  & $F_0$ & 13   & 0.030\% & 0.0011 & 0.606 & 0.46\\
5  & $F_1$ & 117  & 0.267\% & 0.0148 & 0.858 & 0.16\\
5  & $F_2$ & 74   & 0.169\% & 0.0161 & 0.868 & 0.21\\
5  & $F_3$ & 125  & 0.286\% & 0.0288 & 0.883 & 0.10\\
5  & $F_4$ & 450  & 1.088\% & 0.0460 & 0.868 & 0.30\\
\midrule
10 & $F_0$ & 28   & 0.064\% & 0.0014 & 0.483 & 0.46\\
10 & $F_1$ & 228  & 0.520\% & 0.0245 & 0.830 & 0.16\\
10 & $F_2$ & 156  & 0.355\% & 0.0019 & 0.137 & 0.21\\
10 & $F_3$ & 225  & 0.516\% & 0.0285 & 0.868 & 0.10\\
10 & $F_4$ & 827  & 2.087\% & 0.0521 & 0.776 & 0.31\\
\midrule
20 & $F_0$ & 58   & 0.132\% & 0.0011 & 0.395 & 0.46\\
20 & $F_1$ & 448  & 1.022\% & 0.0317 & 0.748 & 0.16\\
20 & $F_2$ & 345  & 0.786\% & 0.0044 & 0.170 & 0.21\\
20 & $F_3$ & 398  & 0.912\% & 0.0460 & 0.815 & 0.10\\
20 & $F_4$ & 1550 & 4.290\% & 0.0806 & 0.713 & 0.33\\
\bottomrule
\end{tabular}
\end{center}
\section{Pooled-score uncertainty}
\label{app:pooled}
\vspace{-6pt}

\begin{center}

\captionof{table}{Moving-block bootstrap uncertainty, concatenated fold-local calibrated scores.}
\label{tab:pooled}
\scriptsize
\begin{tabular}{crrrll}
\toprule
$h$ & Pooled rows & Positives & Prevalence & Pooled AP (95\% CI) & Pooled ROC-AUC (95\% CI)\\
\midrule
5  & 216,723 & 779   & 0.0036 & 0.0049 [0.0029, 0.0082] & 0.541 [0.475, 0.627]\\
10 & 214,978 & 1,464 & 0.0068 & 0.0075 [0.0041, 0.0119] & 0.561 [0.504, 0.655]\\
20 & 211,488 & 2,799 & 0.0132 & 0.0125 [0.0078, 0.0186] & 0.491 [0.454, 0.596]\\
\bottomrule
\end{tabular}
\end{center}

\section{Remaining subset-method grids}
\label{app:grids}
\vspace{-6pt}

\begin{center}
\scriptsize
\captionof{table}{Remaining subset-method grids: subset minus full-data average precision, by horizon
and retained fraction. Uniform and class-stratified are proportional-allocation designs
(Corollary~\ref{cor:balopt}) and $K$-means falls outside the framework
(Proposition~\ref{prop:kmeans}); $K$-means is evaluated at three fractions only. End-to-end speedups
are $1.84$--$2.72\times$ for uniform, $2.03$--$2.89\times$ for class-stratified, and
$0.038$--$0.463\times$ for $K$-means, which is slower than full-data training at every setting.}
\label{tab:uniform}
\begin{tabular}{@{}l rrr rrr rrr@{}}
\toprule
& \multicolumn{3}{c}{$h=5$} & \multicolumn{3}{c}{$h=10$} & \multicolumn{3}{c}{$h=20$}\\
\cmidrule(lr){2-4}\cmidrule(lr){5-7}\cmidrule(lr){8-10}
Retained & Unif. & Strat. & $K$-m. & Unif. & Strat. & $K$-m. & Unif. & Strat. & $K$-m.\\
\midrule
1\%   & $-0.0066$ & $-0.0081$ & $-0.0056$ & $-0.0020$ & $-0.0012$ & $+0.0006$ & $-0.0014$ & $-0.0037$ & $-0.0046$\\
2.5\% & $-0.0068$ & $-0.0046$ & $-0.0070$ & $-0.0014$ & $-0.0037$ & $+0.0003$ & $-0.0049$ & $-0.0036$ & $-0.0039$\\
5\%   & $-0.0036$ & $-0.0061$ & $-0.0084$ & $-0.0001$ & $-0.0025$ & $+0.0028$ & $-0.0026$ & $-0.0014$ & $-0.0026$\\
10\%  & $-0.0076$ & $-0.0045$ & & $+0.0013$ & $+0.0027$ & & $-0.0027$ & $-0.0035$ & \\
20\%  & $-0.0031$ & $-0.0036$ & & $+0.0013$ & $+0.0032$ & & $+0.0025$ & $-0.0003$ & \\
\bottomrule
\end{tabular}
\end{center}

\end{document}